\documentclass[draft, 10pt]{article}
\usepackage{amssymb, amsthm, amsmath}
\numberwithin{equation}{section}

\newcommand{\R}{\mathbb{R}}
\newcommand{\C}{\mathbb{C}}
\newcommand{\Z}{\mathbb{Z}}
\newcommand{\N}{\mathbb{N}}
\newcommand{\Q}{\mathbb{Q}}
\newcommand{\T}{\mathbb{T}}
\newcommand{\Grass}{\mathbb{G}}
\newcommand{\bfP}{\mathbb{P}}
\newcommand{\A}{\mathbb{A}}
\newcommand{\F}{\mathbb{F}}

\newcommand{\e}{\emph}
\newcommand{\rom}{\mathrm}
\newcommand{\ov}{\overline}
\newcommand{\ma}{\mathbf}

\newcommand{\mcal}{\mathcal}
\newcommand{\fr}{\mathfrak}

\newcommand{\ben}{\begin{enumerate}}
\newcommand{\een}{\end{enumerate}}
\newcommand{\eit}{\begin{itemize}}

\newcommand{\ve}{\varepsilon}

\newcommand{\al}{\alpha}
\newcommand{\D}{\Delta}
\newcommand{\del}{\delta}
\newcommand{\om}{\omega}

\newcommand{\be}{\beta}
\newcommand{\la}{\lambda}

\newcommand{\lb}{\left(}
\newcommand{\rb}{\right)}
\newcommand{\lab}{\label}

\newtheorem{thm}{Theorem}
\newtheorem{lem}{Lemma}
\newtheorem*{lemma}{Lemma}
\newtheorem{pro}{Proposition}
\newtheorem{cor}{Corollary}
\newtheorem*{hyp}{Hypothesis}
\newtheorem{con}{Conjecture}

\theoremstyle{definition}
\newtheorem*{ack}{Acknowledgement}

\newcommand{\hcf}{\rom{h.c.f.}}
\renewcommand{\mod}{\hspace{-1mm}\pmod}
\newcommand{\colt}[2]{\genfrac{}{}{0pt}{1}{#1}{#2}}

\newcommand{\x}{\ma{x}}

\newcommand{\y}{\ma{y}}

\newcommand{\bt}{\ma{t}}

\newcommand{\af}{\mathsf{aff}}
\newcommand{\as}{\mathsf{ASH}}
\newcommand{\ho}{\mathfrak{h}}
\newcommand{\conj}{Conjecture }

\begin{document}

\title{Counting Rational Points on Algebraic Varieties}

\author{T.D. Browning$^1$, D.R. Heath-Brown$^2$ and P. Salberger$^3$\\
\small{$^{1,2}$\emph{Mathematical Institute,
24--29 St. Giles',
Oxford OX1 3LB}}\\
\small{$^3${\em Chalmers University of
    Technology, G\"oteborg SE-412 96}}\\
\small{$^1$browning@maths.ox.ac.uk, 
$^2$rhb@maths.ox.ac.uk,  $^3$salberg@math.chalmers.se}}
\date{}

\maketitle

\begin{abstract}

For any $N \geq 2$, let $Z \subset \mathbb{P}^N$ be a geometrically integral
algebraic variety of degree $d$. 
This paper is concerned with the number
$N_Z(B)$ of $\mathbb{Q}$-rational points on $Z$ which have height at most $B$.
For any $\varepsilon>0$ we establish the estimate 
$$
N_Z(B)=O_{d,\varepsilon,N}(B^{\dim Z+\varepsilon}),
$$
provided that $d \geq 6$.  As indicated, 
the implied constant depends at most upon $d, \varepsilon$ and
$N$.\\
Mathematics Subject Classification (2000):  11G35 (14G05) 
\end{abstract}

\section{Introduction}

For any $n \geq 2$, let $F \in \Z[X_0,\ldots,X_n]$ be a form of degree
$d$, that produces an integral hypersurface $F=0$ in $\bfP^{n}$.
Throughout this paper we shall take integrality to mean geometric
integrality,  and likewise we shall say that a polynomial is 
irreducible if it is absolutely irreducible. 
Let $|\x|$ denote the norm $\max_{0\leq i\leq n}|x_i|$ on $\R^{n+1}$, 
and define the quantity 
$$
N(F;B)=\#\{\x \in \Z^{n+1}: F(\x)=0, ~\hcf(x_0,\ldots,x_n)=1, ~|\x| \leq B\},
$$
for any $B\geq 1$.  This paper is motivated by the following
conjecture due to the second author \cite[\conj $2$]{annal}. 

\begin{con}\lab{hb}  Let $\ve>0$ and suppose that $F \in \Z[X_0,\ldots,X_n]$
is an irreducible form of degree $d \geq 2$.
Then we have
$$
N(F;B)=O_{d,\ve,n}(B^{n-1+\ve}).
$$
\end{con}

In all that follows the implied constant in any estimate is absolute unless 
explicitly indicated otherwise.  In the case of \conj \ref{hb}, the
constant is clearly permitted to depend only upon 
$d$, $\ve$ and $n$.

\conj \ref{hb} has already been established  in the cases $n\leq
3$  \cite[Theorems 3 and 9]{annal} and in the case $d=2$ for any $n
\geq 2$  \cite[Theorem 2]{annal}.  Moreover, upon combining recent
work of Broberg and the third author
\cite{3fold1} with work of the first two authors on cubic threefolds
\cite[Theorem 3]{bhb}, the conjecture has also been established in the case $n=4$.
There exist examples for every $n \geq 3$ and $d \geq 2$ for which 
$N(F;B) \gg B^{n-1}$.  To see this it suffices to consider forms of the shape 
$X_0F_0-X_1F_1,$
with $F_0,F_1 \in \Z[X_0,\ldots,X_n]$ of degree $d-1$.
This shows that \conj \ref{hb} is essentially best possible.

It is frequently more useful to have estimates for the corresponding
problem in which $Z \subset \bfP^N$ is an arbitrary integral
algebraic variety. Throughout our work we shall always assume that
$Z$ is a proper subvariety of $\bfP^N$ and that the ideal
of $Z$ is generated by forms defined over $\ov{\Q}$. 
In practice the most interesting situation is
when the variety under consideration is the zero locus of forms which
are all defined over $\Q$, as in the statement
of \conj \ref{hb}.  
Indeed, whenever the ideal of $Z$ is not invariant under the action of
the absolute galois group $\rom{Gal}(\ov{\Q}/\Q)$, the set of rational
points on $Z$ 
automatically lies on a proper subvariety of $Z$.

Let $\x=(x_0,\ldots,x_N) \in \Z^{N+1}$ be any vector such that
$\hcf(x_0,\ldots,x_N)=1$.  
Then we  let $x=[\x]$ denote the corresponding point in $\bfP^N(\Q)$, and write 
$H(x)=|\x|$ for its height.  Conversely, we shall always represent
a point $x\in \bfP^N(\Q)$ by an $(N+1)$-tuple 
$\x=(x_0,\ldots,x_N)\in \Z^{N+1}$  such that
$\hcf(x_0,\ldots,x_N)=1$.  
We proceed by  defining the counting function
\begin{equation}\lab{flu}
N_Z(B)=\#\{x \in Z\cap \bfP^{N}(\Q): ~H(x) \leq B\},
\end{equation}
for any $B \geq 1$ and any algebraic variety $Z \subset \bfP^N$.  
Whenever $X\subset \bfP^n$ is a hypersurface defined by the
equation $F=0$, this definition of $N_X(B)$  is consistent with the
definition of $N(F;B)$.  Indeed, we then have 
$$
N_X(B)=\frac{1}{2}N(F;B),
$$ 
since $\x$ and $-\x$ represent the same
point in projective space.  Note that a hypersurface in $\bfP^{n}$
has dimension $n-1$.  
At this point it is worth highlighting a certain convention concerning 
the letters $n,N$ that we shall follow throughout this paper. 
Whereas $N$ will be used to denote the ambient dimension of a 
projective variety with arbitrary codimension, $n$ will always be used
to denote the ambient dimension of a hypersurface. 
We may now phrase a further conjecture,
which generalises \conj \ref{hb} to arbitrary algebraic varieties.  

\begin{con}\lab{hb'}  Let $\ve>0$ and suppose that $Z \subset \bfP^N$ is an integral
variety of degree $d\geq 2$ and dimension $m$.  
Then we have
$$
N_Z(B)=O_{d,\ve,N}(B^{m+\ve}).
$$
\end{con}

Our first result enables us to relate \conj \ref{hb} and
\conj \ref{hb'}.  By employing 
a straightforward birational projection argument, we shall establish the
following result in \S \ref{pf-proj}.  

\begin{thm}\lab{proj}
Let $\ve>0$.  
Suppose that there exists $\theta_{d,n} \in \R$ such that
$$
N(F;B)=O_{d,\ve,n}(B^{\theta_{d,n}+\ve}),
$$
for any irreducible form $F \in \Z[X_0,\ldots,X_n]$
of degree $d \geq 2$.  Then 
$$
N_Z(B) =O_{d,\ve,N}(B^{\theta_{d,n}+\ve}),   
$$
for any integral variety $Z \subset \bfP^N$ of degree $d \geq 2$
and dimension $n-1$.
\end{thm}

In view  of our remarks above concerning \conj
\ref{hb}, we may deduce from Theorem \ref{proj} that
\conj \ref{hb'} holds for arbitrary projective varieties of
dimension $m \leq 3$, or for arbitrary projective quadrics.
In fact it is clear that Theorem \ref{proj} has the following consequence.

\begin{cor}
\conj \ref{hb} is equivalent to \conj \ref{hb'}.
\end{cor}

For varieties $Z \subset \bfP^N$ of degree $d \geq 3$ and dimension $m
\geq 4$, the best approximation to \conj \ref{hb'} that we have is the estimate
\begin{equation}\lab{pila-1}
N_Z(B) =O_{d,\ve,N}(B^{m+1/d+\ve}).
\end{equation}
This is due to Pila \cite{pila}.  
The primary goal of this paper is to improve upon this
estimate.  In fact our method is largely inspired by the proof of
(\ref{pila-1}), and it will be useful to recall the basic idea here.  By Theorem
\ref{proj} it suffices to consider the case of hypersurfaces defined
by an  irreducible form $F \in \ov{\Q}[X_0,\ldots,X_n]$ of degree
$d\geq 3$.

Now for any $\nu \geq 2$, let $f \in \ov{\Q}[T_1,\ldots,T_\nu]$ be a 
non-zero polynomial of total degree $\delta$, which produces a hypersurface $f=0$ in
$\A^{\nu}$.  We define the counting function 
$$
M(f;B)=\#\{\bt \in \Z^\nu: f(\bt)=0, ~|\bt| \leq B\},
$$
for any $B \geq 1$.
Then the main trick behind the proof of (\ref{pila-1}), as mimicked in 
\S \ref{hypersurfaces} below, is the observation that
\begin{equation}\lab{trick}
N(F;B) \leq \sum_{|b|\leq B} M(f_b;B),
\end{equation}
where $f_b=f_b(T_1,\ldots,T_{n})$ denotes the polynomial
$F(b,T_1,\ldots,T_n)$.  Pila proceeds \cite[Theorem A]{pila} by proving the upper 
bound
\begin{equation}\lab{pilaA}
\#\{\ma{t} \in T\cap \Z^\nu: ~|\bt| \leq B\}
=O_{\delta,\ve,\nu}(B^{\mu-1+1/\delta+\ve}),
\end{equation}
for any integral affine variety $T \subset \A^{\nu}$ 
of degree $\delta$ and dimension $\mu$. 
Such an estimate yields the corresponding bound $O_{d,\ve,n}(B^{n-2+1/d+\ve})$
for $M(f_b;B)$, provided that $f_b$ is irreducible and has degree
$d$.  Upon treating the remaining degenerate cases separately, this
suffices to establish (\ref{pila-1}).

The proof of (\ref{pilaA}) is argued by induction on the dimension
$\mu$ of $T$.  Pila starts the induction at  $\mu=1$, for which he
employs his earlier joint work with Bombieri \cite{bp}.  
Our ability to improve upon (\ref{pila-1}) ultimately stems from the 
idea of using a similar induction argument to that of Pila, but
instead taking the case $\mu=2$ as the inductive base.  
In order to simplify the exposition of our work, it will be convenient to make the
following hypothesis.

\begin{hyp}[Affine surface hypothesis]
Let $\ve>0$ and let $\mcal{P}$ be any set of polynomials 
$f \in \Z[T_1,T_2,T_3]$ of degree $\delta$.
Then there exists $\al \in \R$ such that
$$
M(f;B)=O_{\delta,\ve}(B^{\al+\ve}),
$$
for any $f \in \mcal{P}$.
\end{hyp}

As highlighted above, the implied constant here is only permitted to depend upon
$\del$ and $\ve$, and not on the individual coefficients of $f$.
We shall henceforth write $\as[\al, \mcal{P}]$ to denote the affine surface
hypothesis holding with exponent $\al$ for the set of polynomials $\mcal{P}$.
Let $\mcal{I}_\delta$ denote the set of irreducible polynomials in 
$\Z[T_1,T_2,T_3]$ of degree $\delta$.  Then we have already seen  via Pila's estimate (\ref{pilaA}) that
$\as[1+1/\del, \mcal{I}_\delta]$ holds. In fact it is easy to construct examples which show that this result
is essentially best possible.  Thus the irreducible
polynomial $T_1-T_2^\del$ has $\gg B^{1+1/\del}$ zeros $\bt\in \Z^3$
with $|\bt|\leq B$.  However it is possible to take smaller
exponents in the affine surface hypothesis by restricting attention to
a suitable class of polynomials.  
This basic idea will enable us to obtain compelling evidence for \conj \ref{hb'}.
For any polynomial $g \in \Z[T_1,\ldots,T_\nu]$, let
$\ho(g)$ denote the homogeneous part of $g$ of
maximal degree.  We let  $\mcal{J}_\delta$ denote the set of
polynomials $f \in \Z[T_1,T_2,T_3]$ of degree $\delta$, such that
$\ho(f)$ is irreducible.  In particular it is clear that
$\mcal{J}_\delta\subset \mcal{I}_\delta$.
With this in mind, the following result will
be established in \S \ref{aff-surface}.

\begin{thm}\lab{ash-2}
Let $\del \geq 3$. Then $\as[\al, \mcal{J}_\del]$ holds with 
$$
\al=  \left\{
\begin{array}{ll}
5/(3\sqrt{3})+1/4, & \del=3,\\
3/(2\sqrt{\del})+1/3, & \mbox{$\del=4$ or $5$},\\
1, & \del\geq 6.
\end{array}
\right.
$$ 

\end{thm}

We note for comparison that $\ho(T_1-T_2^\del)$ is a product of $\del$
linear factors, so that $T_1-T_2^\del \not\in \mcal{J}_\del$.
In view of Theorem \ref{ash-2}, it seems reasonable to formulate the following conjecture.

\begin{con}\lab{con-ash}
Let $\del \geq 2.$  Then $\as[1,\mcal{J}_\del]$ holds.
\end{con}

We now examine how Conjecture \ref{con-ash} relates to Conjecture
\ref{hb'}.  The following key result will be established in \S \ref{hypersurfaces}.

\begin{thm}\lab{imply-2}
Let $\ve>0$ and suppose that $Z \subset \bfP^N$ is an integral 
variety of degree $d\geq 2$ and dimension $m$.  
Then if $\as[\al,\mcal{J}_d]$ holds we have
$$
N_Z(B)=O_{d,\ve,N}(B^{m-1+\al+\ve}).
$$
In particular \conj \ref{con-ash} implies \conj \ref{hb'}.
\end{thm}

As a direct consequence of Theorem \ref{imply-2}, we may clearly apply
Theorem \ref{ash-2} to deduce the following estimate.

\begin{cor}\lab{cor2}
Let $\ve>0$ and suppose that $Z \subset \bfP^N$ is an integral 
variety of degree $d\geq 3$ and dimension $m$.  
Then we have
$$
N_Z(B)\ll_{d,\ve,N}\left\{
\begin{array}{ll}
B^{m-3/4+5/(3\sqrt{3})+\ve}, & d=3,\\
B^{m-2/3+3/(2\sqrt{d})+\ve}, &  \mbox{$d=4$ or $5$},\\
B^{m+\ve}, & d\geq 6.
\end{array}
\right.
$$
\end{cor}

It is worthwhile highlighting that Corollary \ref{cor2} improves upon 
(\ref{pila-1}) for all values of $d\geq 3$, and plainly establishes
Conjecture \ref{hb'} for $d\geq 6$.   We have already seen that
\conj \ref{hb'} holds for $d=2$.

We end this introduction by providing a brief summary of the contents
of this paper.
The following section will contain a review of some basic properties of 
Grassmannians, in addition to 
presenting a number of background estimates that will be useful to us.
In \S \ref{pf-proj} a birational projection argument will be used to
establish Theorem \ref{proj}, and in \S \ref{hypersurfaces} we shall
adapt Pila's inductive proof of 
(\ref{pila-1}) to obtain a proof of Theorem \ref{imply-2}.
The remaining sections will all be taken up with establishing Theorem \ref{ash-2}.
Thus in \S \ref{geometry} we shall collect together some preliminary results from
algebraic geometry, before continuing with the
proof of Theorem \ref{ash-2} proper in \S \ref{aff-surface}.

\begin{ack}
The authors are very grateful to the referees, for their careful reading of the text and numerous pertinent
suggestions. While working on this paper, the first
author was supported by
EPSRC grant number GR/R93155/01.
\end{ack}

\section{Preliminary estimates}\lab{prelim}

We begin by establishing  the following
simple estimate.   This will be crucial to the proof of Theorem
\ref{ash-2} in \S \ref{aff-surface} below.

\begin{lem}\lab{4count}
Let $p(t)=a_0+a_1t+\cdots+a_\delta t^\delta \in \C[t]$ be a polynomial
of degree $\delta$.  Then for any $T \geq 1$ we have
$$
\#\{t \in \Z: |p(t)| \leq T\} \ll_\delta 1+ (T/|a_\delta|)^{1/\delta}.
$$ 
\end{lem}
 
\begin{proof}
Our starting point in the proof of Lemma \ref{4count} is the
trivial  estimate
$$
\#\{t \in \Z: |t-\la| \leq L\} \ll 1+L,
$$
for any $\la \in \C$ and $L\geq 0$.  Here the implied
constant is independent of $\la$.
Next we write $p(t)=a_\delta q(t)$, say, with $q \in \C[t]$ a monic polynomial.  
Then we may factorise $q$ as a product of linear polynomials
$$
q(t)=\prod_{i=1}^\delta (t-\la_i),
$$ 
for $\la_1,\ldots, \la_\delta \in \C$.
Let $S(T;p)$ denote
the set of $t \in \Z$ for which $|p(t)| \leq T$.  Then
we clearly have $S(T;p)=S(T/|a_\delta|;q)$.  We proceed by sorting the points of
$S(T/|a_\delta|;q)$ into sets $S_1,\ldots,S_\delta$, according to the value
of $i\in \{1,\ldots,\delta\}$ for which $|t-\la_i|$ is least.
But then it follows that for each $1 \leq i \leq \delta$ we have
$$
\#S_i \leq \{t \in \Z: |t-\la_i| \leq (T/|a_\delta|)^{1/\delta}\} \ll 1+ (T/|a_\delta|)^{1/\delta},
$$
by the above.  Since $S(T/|a_\delta|;q)=\bigcup_{i=1}^\delta S_i$, the
result follows.
\end{proof}

Recall the notation (\ref{flu}) for the counting function attached to
a variety $Z\subset \bfP^N$.  
The following ``trivial'' 
estimate is established in \cite[Theorem 1]{bhb}.

\begin{lem}\lab{triv}
Let $Z \subset \bfP^{N}$ be a variety of degree $d$ and
dimension $m$.   Then we have
$$
N_Z(B) =O_{d,N}(B^{m+1}).
$$
\end{lem}

It is easy to see that Lemma \ref{triv} is best possible whenever $Z$
contains a linear subspace of dimension $m$ that is defined over
$\Q$.

A useful consequence of Lemma \ref{triv} is that for any 
variety $Z \subset \bfP^{N}$ of degree $d$ and
dimension $m<N$, one can always find a point $x \in \bfP^N(\Q)$ such
that $x \not\in Z$ and $H(x)=O_{d,N}(1)$.  
We now discuss how this fact can
be generalised somewhat.
For non-negative integers $k<N$, let $\Grass(k,N)$ denote the Grassmannian
which parameterises $k$-planes contained in $\bfP^{N}$. 
Here a $k$-plane is merely a linear subspace of $\bfP^N$ of dimension
$k$. It is well-known that $\Grass(k,N)$ can be embedded in $\bfP^\nu$ via
the Pl\"ucker embedding, where
$$
\nu=\binom{N+1}{k+1}-1.
$$
To see this explicitly, we note that since any $k$-plane $M \subset \bfP^N$ is spanned by $k+1$
points, we may represent $M$ as a $(k+1)\times (N+1)$ matrix $\ma M$
of rank $k+1$.  The Pl\"ucker embedding is then the map which
takes this matrix representation $\ma M$ to the point $[\det
\ma{M}_0, \ldots, \det \ma{M}_{\nu}] \in \bfP^{\nu}$, where each $\ma{M}_i$ is the
square matrix composed of any $k+1$ columns of $\ma M$.   In fact
this definition is only valid up to signs, but it is satisfactory for our
purposes.  One then has that $\Grass(k,N)$ is a 
subvariety of $\bfP^\nu$ of  dimension 
$(k+1)(N-k).$  With these facts in mind we have the following result.

\begin{lem}\lab{triv-grass}
For non-negative integers $k<N$, let $Y \subset \Grass(k,N)$ be a
proper subvariety of degree $d$.  Then there exists a point $y \in
\Grass(k,N)\cap\bfP^\nu(\Q)$ such that $y \not\in Y$ and $H(y)=O_{d,N}(1)$.
\end{lem}

\begin{proof}
Whenever $k=0$ this reduces to the statement that we can
always find points of small height in
$\bfP^N(\Q)$ which lie off a given proper subvariety.
This much has already been observed as a consequence of Lemma
\ref{triv}.

For the general case $k \geq 1$, we write $\phi: (\bfP^N)^{k+1}
\rightarrow \bfP^\nu$ for the rational map taking $(k+1)$-tuples of
points in $\bfP^N$ to the point in $\bfP^{\nu}$ formed from the
maximal minors of the corresponding $(k+1)\times (N+1)$ matrix.
If $U \subset (\bfP^N)^{k+1}$ is the open subset of all $(k+1)$-tuples
which span a $k$-plane, then $\phi:U \rightarrow \bfP^\nu$ is a morphism
whose image $\phi(U)$ is the Grassmannian $\Grass(k,N) \subset \bfP^\nu$.
We are interested in the inverse image $U_Y=\phi^{-1}(Y)$ in $U$.  
In particular it is clear that there exist non-zero multi-homogeneous 
polynomials  
$$
F_i(X_0^{(0)},\ldots,X_N^{(0)};\ldots;X_0^{(k)},\ldots,X_N^{(k)}), \quad (1\leq i \leq t),
$$
say, which generate the ideal of $U_Y$.

Let $\x^{(\ell)}=(X_0^{(\ell)},\ldots,X_N^{(\ell)})$ for $0\leq \ell \leq k$.  
Then for any $1 \leq i \leq t$ we may write
$$
F_i(\x^{(0)};\ldots;\x^{(k)})=\sum_j F_{ij}(\x^{(0)})G_{ij}(\x^{(1)};\ldots;\x^{(k)}),
$$
for appropriate forms $F_{ij}(\x^{(0)})$ not all identically zero, and
linearly independent polynomials $G_{ij}(\x^{(1)};\ldots;\x^{(k)})$.
In view of Lemma \ref{triv} we may choose a point
$\ma{a}^{(0)}\in\Z^{N+1}$ such that $F_{ij}(\ma{a}^{(0)})$ are not all
zero and
$|\ma{a}^{(0)}|=O_{d,N}(1)$.  It follows that the polynomial
$F_i(\ma{a}^{(0)};\ma{x}^{(1)};\ldots;\x^{(k)})$ does not vanish identically.
Thus we may proceed by writing 
$$
F_i(\ma{a}^{(0)};\ma{x}^{(1)};\ldots;\x^{(k)})=\sum_j
\tilde{F}_{ij}(\ma{a}^{(0)};\x^{(1)})\tilde{G}_{ij}(\x^{(2)};\ldots;\x^{(k)}), 
$$
with $\tilde{F}_{ij}(\ma{a}^{(0)};\x^{(1)})$ not all identically
zero and $\tilde{G}_{ij}(\x^{(2)};\ldots;\x^{(k)})$ linearly independent.  We then choose a point
$\ma{a}^{(1)}\in\Z^{N+1}$ such that $\ma{a}^{(0)}$ and $\ma{a}^{(1)}$ are linearly independent, the numbers 
$F_{ij}(\ma{a}^{(0)};\ma{a}^{(1)})$ are not all zero, and $|\ma{a}^{(1)}|=O_{d,N}(1)$.
Lemma \ref{triv} ensures that this is possible.  Continuing in this
fashion we ultimately obtain a $(k+1)$-tuple of rational points $(a^{(0)},\ldots,a^{(k)}) \in
(\bfP^N)^{k+1}\setminus U_Y$, which span a $k$-plane and have heights
$$
H(a^{(0)}),\ldots,H(a^{(k)})=O_{d,N}(1).
$$ 
The point $y=\phi(a^{(0)},\ldots,a^{(k)}) \in
\Grass(k,N)$ is sufficient for Lemma \ref{triv-grass}.
\end{proof}

We end this section by recalling two further estimates 
that will be useful in \S \ref{aff-surface}.  The first of these is a rather general estimate
for the number of rational points of bounded height on a projective
plane curve \cite[Theorem 3]{annal}.

\begin{lem}\lab{bub}
Let $\ve>0$ and suppose that $C \subset \bfP^2$ is an integral 
curve of degree $d$.
Then we have
$$
N_C(B)=O_{d,\ve}(B^{2/d+\ve}).
$$
\end{lem}

Let $H(G)$ denote the maximum modulus of the coefficients of
a form $G \in \Z[X_0,\ldots,X_n]$, and say that $G$ is primitive if the
highest common factor of its coefficients is $1$.  The following
result is established by exactly the same argument as for \cite[Theorem 4]{annal}.

\begin{lem}\lab{theta}
Let $G \in \Z[X_0,\ldots,X_n]$ be a primitive non-zero form of degree
$d \geq 2$, defining a hypersurface $Z\subset \bfP^n$.  Then either
$$
H(G) \ll B^{\theta}, \quad \theta= d\left(\colt{d + n}{n}\right),
$$
or else  there exists a form
$G' \in \Z[X_0,\ldots,X_n]$ of degree $d$, not proportional to $G$,
such that $G'$ vanishes at each point of the set $\{x\in
Z\cap\bfP^n(\Q): H(x)\leq B\}$.
\end{lem}

\section{Varieties in $\bfP^N$}\lab{pf-proj}

In this section we shall establish Theorem \ref{proj}.  
Let $Z \subset \bfP^{N}$ be an integral variety of degree $d \geq
2$ and dimension $m$.    
The thrust of our work will be taken up with 
constructing a projection $Z \rightarrow \bfP^{m+1}$, 
such that the image $\ov{Z} \subset \bfP^{m+1}$ 
is an  integral variety of degree $d$ and dimension $m$.  Moreover
we shall want to choose our projection in such a way that we have
\begin{equation}\lab{aim}
N_Z(B) \leq d N_{\ov{Z}}(cB),
\end{equation}
for some $c \ll_{d,N} 1$.
Once this is accomplished, the statement of Theorem \ref{proj} easily
follows since now $\ov{Z} \subset \bfP^{m+1}$ is an integral
hypersurface of degree $d$.
We claim that such a projection always exists.

We must first deal with the possibility that $Z$ is ``degenerate'', by
which we shall mean that $Z \subset H$ for some hyperplane $H \in \Grass(N-1,N)$.
Assume that $H$ is defined by the
linear equation $\sum_{i=0}^N a_i X_i=0$, and suppose 
without loss of generality that $a_{N}\neq 0$.  
Then the point $y=[0,\ldots,0,1]$ is not
contained in $H$ and so the projection $\pi_y$ from $Z$ onto the
hyperplane $X_N=0$ is a regular map.  Moreover $\pi_y$ is clearly
birational onto the image $\ov{Z}=\pi_y(Z) \subset \bfP^{N-1}$, and so $\ov{Z}$ is an
integral variety of degree $d$ such that  $N_Z(B) \leq N_{\ov{Z}}(B)$.
Now either $\ov{Z}$ is not degenerate, or we may
repeat the argument once again.  
Arguing in this way it suffices to assume henceforth  that $Z$ is
not degenerate.

We shall need to find an $(N-m-2)$-plane $\Lambda \subset
\bfP^N$ such that the projection $\pi_\Lambda: Z \rightarrow \Gamma$ from
$\Lambda$ onto any $(m+1)$-plane $\Gamma$ disjoint from $\Lambda$, is birational onto the image.
Moreover, we shall want $\Lambda$ to be 
defined over $\Q$ and have height $H(\Lambda) \ll_{d,N} 1$.
Here we define the height of $\Lambda$
to be the standard multiplicative height of its coordinates in
$\Grass(N-m-2,N)$, under the Pl\"ucker embedding.  
Now it is well-known that for a generic $\Lambda \in \Grass(N-m-2,N)$
the projection described above is birational.  
Hence the main work comes from showing that one can
always find such a linear space which is defined over $\Q$ and which has
small height.  To do so we shall need to establish the following
preliminary result.

\begin{lem}\lab{dim1}
Let $Z\subset \bfP^N$ be an integral 
variety of degree $d$ and dimension $m$.
Let
$$
Y=\{\Lambda \in \Grass(N-m-2,N): \dim S_{\Lambda,Z} \geq m\}, 
$$
where $S_{\Lambda,Z}$ denotes the set of
$M \in \Grass(N-m-1,N)$ such that $\Lambda \subset M$ and $M$ intersects
$Z$ in at least two (possibly coincident) points.
Then $Y$ is a proper subvariety of $\Grass(N-m-2,N)$ of degree $O_{d,N}(1)$.
\end{lem}

\begin{proof}
It is well known that varieties in $\bfP^N$ of degree $d$ and
dimension $m$ can be parameterised by a quasi-projective variety
$\mcal{C}_{d,m}(\bfP^N)\subset \bfP U$, where $U$ is the vector space of
multi-homogeneous polynomials of multi-degree $(d,\ldots,d)$ in $m+1$
sets of $N+1$ variables.  This is the open Chow variety (see
\cite[Theorem 21.2]{harris}, for example) of varieties of degree
$d$ and dimension $m$ in $\bfP^N$. If $V\subset \bfP^N$ is any
variety of degree $d$ and
dimension $m$, we let $[c_V]$ denote the corresponding Chow point in
$\mcal{C}_{d,m}(\bfP^N)$.  
We shall also denote by $F_{d,m}$ the universal family of all
$m$-dimensional subvarieties of degree $d$ in $\bfP^N$.  It is the
closed subvariety $F_{d,m}\subset \mcal{C}_{d,m}(\bfP^N)\times \bfP^N$
of points $([c_V],x)\in \mcal{C}_{d,m}(\bfP^N)\times \bfP^N$ such that
$x\in V$.
Much as we have followed the convention of
identifying linear spaces in $\bfP^N$ with the corresponding Pl\"ucker
point in the Grassmannian, we shall
frequently write $V$ to denote the Chow point $[c_V]$.

Write $\mcal{C}=\mcal{C}_{d,m}(\bfP^N)$,  $G_1=\Grass(N-m-2,N)$ and
$G_2=\Grass(N-m-1,N)$  for convenience.
We proceed by considering the incidence correspondence
$$
I= \left\{(\Lambda,M,V) \in G_1\times G_2 \times \mcal{C}: \Lambda
  \subset M, ~\#(M\cap V)\geq 2 \right\}.
$$
Here the condition $\#(M\cap V)\geq 2$ asserts that $M$ intersects
$V$ in at least two (possibly coincident) points. We formally
define $\#(M\cap V)$ to be infinity if the intersection of $M$ and $V$
contains a component of dimension $\geq 1$.
Let $I_1\subset G_1\times G_2$ be the set of 
$(\Lambda,M) \in G_1\times G_2$ such that $\Lambda \subset M$, and let
$I_2\subset G_2\times \mcal{C}$ be the set of 
$(M,V) \in G_2 \times \mcal{C}$ such that $\#(M\cap V)\geq 2$.  Then
if $\phi_1$ denotes the projection  $G_1\times G_2 \times
\mcal{C}\rightarrow  G_1\times G_2$, and 
$\phi_2$ the projection  $G_1\times G_2 \times
\mcal{C}\rightarrow  G_2\times \mcal{C}$, then
it is plain that $I=\phi_1^{-1}(I_1)\cap
\phi_2^{-1}(I_2)$.
Thus in order to show that $I$ is closed in the Zariski topology, it
suffices to show that $I_1$ and $I_2$ are both closed sets.  
It is well-known that $I_1$ is closed (see \cite[Proposition
I.9.9.3]{ega1}, for example). 
To see that $I_2$ is closed in $G_2\times
\mcal{C}$ we show that the intersection cardinality of $M$ and $V$ is an upper semicontinuous
function.  Let $\Omega_1$ be the set of $(M,V)\in G_2\times\mcal{C}$ such
that $M\cap V$ contains a component of dimension $\geq 1$, and let
$\Omega_2=I_2\setminus \Omega_1$.
It is clear that $\Omega_1$ is closed, and so it remains to show that
$\#(M\cap V)\geq n$ gives a closed condition on $\Omega_2$, for any
positive integer $n$.  To do so we apply
\cite[Excercise II.5.8(a)]{hart} to the restriction
of the coherent sheaf $\phi_*\mcal{O}_Z$ to $\Omega_2$ for the projection
$\phi:Z\rightarrow G_2\times \mcal{C}$ from the closed subset $Z\subseteq
G_2\times F_{d,n}$ of all $(M, V,x)\in G_2\times F_{d,n}$ such that $x\in M$.

Let $\pi: I\rightarrow G_1 \times \mcal{C}$ be the projection
map, and let $W\subset G_1 \times \mcal{C}$ denote the subset over
which the fibre of $\pi$ has dimension at least $m$.
Then $W$ is a closed subvariety of $G_1
\times \mcal{C}$, by the upper semicontinuity of fibre dimension
\cite[Corollaire 13.1.5]{ega4}. Let $[c_Z] \in \mcal{C}$ be
the Chow point corresponding to the variety $Z$.   Then if
$\rho:W \rightarrow \mcal{C}$ denotes the projection onto the
second factor, we consider the fibre
$$
\rho^{-1}([c_Z])=\{\Lambda \in G_1: (\Lambda,Z) \in W\}.
$$
It follows that $\rho^{-1}([c_Z])$ is a projective subvariety of
$G_1=\Grass(N-m-2,N)$, whose degree can be bounded in terms of $d$ and
$N$ alone.  To see this we recall the embedding $G_1\times
\mcal{C} \subset \bfP^r\times \bfP U$, with $r=\binom{N+1}{N-m-1}-1$,
so that $W\subset \bfP^r\times
\bfP U$ is cut out by $O_{d,N}(1)$ multi-homogeneous polynomials of multi-degree
$O_{d,N}(1)$.  The claim then follows on specialising at $[c_Z]$, and
applying B\'ezout's theorem (in the form given by \cite[Theorem
8.4.6]{fulton}). Finally we note that 
$Y=\rho^{-1}([c_Z])$ in the statement of the lemma. 

In order to complete the proof of Lemma \ref{dim1} we need to show
that $Y$ is a proper subvariety of $G_1$.  As in the statement of
the lemma let $S_{\Lambda,Z}$ denote the set of  $(N-m-1)$-planes containing $\Lambda$
and meeting $Z$ in at least two points, for any $\Lambda \in G_1$.  
Now let $\Lambda\in G_1$ be such that $\Lambda \cap Z$ is empty,
and let $\rho_\Lambda: Z\rightarrow G_2$ be the map given by sending a
point $x \in Z$ to the projective linear subspace $\langle
x,\Lambda\rangle$ spanned by $x$ and $\Lambda$.
Then $\rho_\Lambda$ is birational onto its image if and only if $\dim
S_{\Lambda,Z} <m$.   This follows on combining the observation that
$S_{\Lambda,Z}\subset \rho_\Lambda(Z)$ is the proper closed subset of $M\in \rho_\Lambda(Z)$ such that
$\rho_\Lambda^{-1}(M)$ consists of more than one point, with the fact
that $\rho_\Lambda(Z)$ is integral and has dimension $m$.
Now let $\Lambda \in G_1$ be generic.  Then it is well known that the map
$\rho_\Lambda$ is birational onto its image \cite[Exercise
11.23]{harris}.  We therefore conclude that $\Lambda \not\in Y$, which 
thereby completes the proof of Lemma~\ref{dim1}.
\end{proof}

We now proceed with the proof of Theorem \ref{proj}.  
Let $Y'$ denote the set of $\Lambda \in
\Grass(N-m-2,N)$ which meet $Z$.  Then $Y'$ is a proper subvariety of
$\Grass(N-m-2,N)$ and has degree $O_{d,N}(1)$.
In view of Lemma \ref{triv-grass}, we may therefore conclude from Lemma \ref{dim1} that there exists 
an $(N-m-2)$-plane $\Lambda \subset \bfP^N$ such that
$\Lambda$ is defined over $\Q$ and has height $H(\Lambda) \ll_{d,N} 1$,
and such that $\Lambda \not\in Y \cup Y'$.  In particular $\Lambda
\cap Z$ is empty and the 
projection $\rho_\Lambda$ from $\Lambda$ is birational onto its 
image.  Below we shall select a certain linear space $\Gamma \in \Grass(m+1,N)$
which does not meet $\Lambda$, and which is defined over $\Q$ and has height
$H(\Gamma) \ll_{d,N} 1$.  
For the moment, let $\Gamma \in \Grass(m+1,N)$
be any $(m+1)$-plane that does not meet $\Lambda$, and let
$$
\pi_\Lambda: Z \rightarrow \Gamma
$$
be the projection of $Z$ from $\Lambda$ onto $\Gamma$.  
Then by construction $\pi_\Lambda$ is a regular map that is birational
onto its image, and we may therefore 
conclude that the image $\ov{Z}=\pi_\Lambda(Z) \subset \Gamma$
is an integral variety of dimension $m$ and degree $d$.  
Finally we note that the fibre over any point of $\ov{Z}$ contains at
most $d$ points.

Now suppose that $\Lambda$ is equal to the linear span of $N-m-1$ points
$$
[\ma{h}_{1}],\ldots,[\ma{h}_{N-m-1}] \in \bfP^N,
$$ for vectors $\ma{h}_{1},\ldots,\ma{h}_{N-m-1} \in \Z^{N+1}$ of modulus $O_{d,N}(1)$.
Then by Lemma \ref{triv-grass} we may select a vector $\ma{g}_1 \in
\Z^{N+1}$ of modulus $O_{d,N}(1)$ such that 
$\ma{g}_1.\ma{h}_j = 0$ for $2 \leq j \leq N-m-1$, and
$\ma{g}_1.\ma{h}_1 \neq 0.$   Similarly, for $2 \leq i \leq N-m-1$ we
may find vectors $\ma{g}_i \in
\Z^{N+1}$ of modulus $O_{d,N}(1)$ such that 
$\ma{g}_i.\ma{h}_j = 0$ for $i \neq j$ and $\ma{g}_i.\ma{h}_i \neq
0$.  These vectors are clearly linearly independent, and define the 
$(m+1)$-plane 
$$
\Gamma: \quad \ma{g}_1.\x=\cdots =\ma{g}_{N-m-1}.\x=0.
$$
By construction $\Lambda \cap \Gamma$ is empty.
We may now explicitly write our projection in the form 
$$
\pi_\Lambda(x)=\left[\ma{x} -
  \frac{(\ma{g}_1.\x)}{(\ma{g}_1.\ma{h}_1)}\ma{h}_1 - \cdots -
  \frac{(\ma{g}_{N-m-1}.\x)}{(\ma{g}_{N-m-1}.\ma{h}_{N-m-1})}\ma{h}_{N-m-1} \right], 
$$ 
for any $x=[\x] \in Z$.
Let $\la=\prod_i \ma{g}_{i}.\ma{h}_i$ and let
$\la_j=\prod_{i\neq j }\ma{g}_{i}.\ma{h}_i$.
Then since $|\ma{g}_i.\x| \leq (N+1)|\ma{g}_i||\x|$ for $1 \leq i \leq
N-m-1$, we deduce that $H(\pi_\Lambda(x)) \leq cH(x)$ with
$$
c=|\la| +(N+1)\sum_{i=1}^{N-m-1} |\la_i|H([\ma{g}_i])H([\ma{h}_i])
\ll_{d,N} 1.
$$
It therefore follows that (\ref{aim}) holds, as claimed, with $\ov{Z}=\pi_\Lambda(Z)
\subset \bfP^{m+1}$.    
This completes the proof of Theorem \ref{proj}.

\section{Affine hypersurfaces}\lab{hypersurfaces}

In this section we shall establish
Theorem \ref{imply-2}.   The main step in the  proof
will involve using the affine surface hypothesis to deduce a corresponding bound for affine hypersurfaces.
In view of Theorem \ref{proj} it will suffice to establish the
estimate
$$
N(F;B)=O_{d,\ve,n}(B^{n-2+\al +\ve}),
$$
for any irreducible form $F \in \Z[X_0,\ldots,X_n]$
of degree $d \geq 2$, provided that $\as[\al,\mcal{J}_d]$ holds.  
Let $X \subset \bfP^{n}$ denote the hypersurface defined by the
equation $F=0$.  
We shall need some control over the hyperplanes $H \in
\Grass(n-1,n)={\bfP^n}^*$ 
for which $X \cap H$ is not integral.  For this we shall
apply the following well-known result (see \cite[Lemma 2.2.1]{3fold1}, for example).

\begin{lem}\lab{elim}
Let $n \geq 3$ and let $X \subset \bfP^{n}$ be an integral hypersurface of
degree $d$.  Then there exists a non-zero form $E \in
\Z[X_0,\ldots,X_n]$ of degree $O_{d,n}(1)$, such that $E(\ma{a})=0$
whenever the hyperplane 
$$
H_{\ma{a}}: \quad a_0X_0+\cdots+a_nX_n=0
$$
produces a hyperplane section $X \cap H_\ma{a}$ that is not integral.
\end{lem}

Lemma \ref{elim} may be viewed as an explicit statement of the basic
geometric fact that the generic hyperplane section of an integral variety
is integral.  
Using Lemma \ref{triv} and Lemma \ref{elim} we may find a vector
$\ma{a} \in \Z^{n+1}$ such that $|\ma{a}| \ll_{d,n} 1$ and $X \cap H_\ma{a}$ is
integral.   After a suitable change of coordinates we
may assume henceforth  that $H_\ma{a}$ is the hyperplane $X_0=0$, and
hence that the form  $F(0,X_1,\ldots,X_n)$ is
irreducible.

For any $b \in \Z$ we define the polynomial
$
f_b(T_1,\ldots,T_{n})=F(b,T_1,\ldots,T_n).
$
Then $f_b \in \Z[T_1,\ldots,T_{n}]$ is a polynomial of degree
$d$ such that $\ho(f_b)$ is irreducible, 
so that $f_b \in \mcal{J}_d$. 
Moreover  (\ref{trick}) implies that
$$
N(F;B) \leq \sum_{|b|\leq B} M(f_b;B).
$$
In order to complete the proof of Theorem \ref{imply-2} it is
therefore sufficient to establish the following result.

\begin{lem}\lab{pilaA-2}
Let $\ve>0$ and suppose that $\as[\al,\mcal{J}_\del]$ holds.  
For any $\nu \geq 3$ let $f \in \Z[T_1,\ldots,T_\nu]$ 
be a polynomial of degree $\del$ such that $\ho(f)$ is irreducible.  
Then we have
$$
M(f;B) \ll_{\del,\ve,\nu} B^{\nu-3+\al+\ve}.
$$
\end{lem}

\begin{proof}
The proof of Lemma  \ref{pilaA-2} will involve taking repeated
hyperplane sections of the affine hypersurface $f=0$.  
We shall argue by 
induction on $\nu$.  The case $\nu=3$ is satisfactory by the
assumption that $\as[\al,\mcal{J}_\del]$ holds. 
For general $\nu\geq 4$ we write $f_0=\ho(f)$, so that
$f_0(T_1,\ldots,T_\nu)$ is an irreducible form of degree $\del$ by 
assumption.  Employing Lemma \ref{triv} and Lemma \ref{elim} we may
therefore assume, possibly after a change of variables, that the
hyperplane section $f_0=T_\nu=0$ produces an integral
hypersurface.  
That is, we may assume that the homogeneous part of the polynomial
$f_c=f(T_1,\ldots,T_{\nu-1},c)$ is irreducible for any $c \in \Z$.
Thus we obtain the bound
$$
M(f;B)\leq \sum_{|c| \leq B} M(f_c; B) \ll_{\del,\ve,\nu} B^{\nu-3+\al+\ve},
$$
by the induction hypothesis.  This completes the proof of Lemma \ref{pilaA-2}.
\end{proof}

\section{Geometry of surfaces}\lab{geometry}

The remainder of the paper will be taken up with the proof of Theorem
\ref{ash-2}.  In this section
we collect together some of the geometric results that
will be necessary to us, throughout which $X \subset \bfP^3$ will
always denote a surface of degree $d \geq 3$.  
Let $H \in {\bfP^3}^*$ be a plane such that $X\cap H$ is integral, 
where ${\bfP^3}^*=\Grass(2,3)$ as above, and 
let $C\subset X$ be an integral curve of degree $e<d$.
In particular it follows that $X$ is integral and that $C$ is not
contained in the plane section $X\cap H$.  For given $y \in X\cap H$,
we shall begin by investigating
the integral curves $C\subset X$ of degree $e<d$ that have
intersection multiplicity $e$ with $H$ at $y$.

\begin{lem}\lab{9cup10}
Let $e<d$ and let 
$H \in {\bfP^3}^*$ be a plane such that $X\cap H$ is integral.
Then for any $y \in X\cap H$ there are $O_{d}(1)$ integral curves
$C\subset X$ of degree $e$ such that $C$ has 
intersection multiplicity $e$ with $H$ at $y$.
\end{lem}

\begin{proof}
Let $\mcal{P}$ be the set of all Hilbert polynomials of integral space
curves of degree $e$.  Then Castelnuovo's inequality implies that
$\mcal{P}$ is finite.  Let $\mcal{H}$ be the Hilbert scheme of all
closed subschemes of $\bfP^3$ with Hilbert polynomials in $\mcal{P}$,
and let $\mcal{H}'$ be the Hilbert scheme of all
closed surfaces $X\subset \bfP^3$ of degree $d$.  
We shall also need to introduce the universal family $\mcal{F}\subset
\bfP^3\times\mcal{H}$ of all closed subschemes of $\bfP^3$ with
Hilbert polynomials in $\mcal{P}$,  and the open subscheme
$\mcal{H}_0$  of $\mcal{H}$ representing curves not in $H$.

For any pair $P=(C,y) \in \mcal{H}_0\times H$, let $i_y(C\cap H)$ denote the intersection
multiplicity of $C$ and $H$ at $y$. Also, let $\mcal{G}\subset \bfP^3
\times \mcal{H} \times H$ be the scheme-theoretic intersection of
$\mcal{F}\times H\subset \bfP^3\times \mcal{H}\times H$ with the
inverse image of the diagonal $H^\Delta \subset \bfP^3 \times H$ under
the projection $\bfP^3\times \mcal{H}\times H\rightarrow \bfP^3\times
H$. Then $i_y(C\cap H)$ can be interpreted as
$\dim_{k(P)}\phi_*\mcal{O}_{\mcal{G}}\times k(P)$  
for the projection $\phi :\mcal{G}\rightarrow \mcal{H}\times H$.
Hence $i_y(C\cap H)$ is an upper semicontinuous function on $\mcal{H}_0
\times H$ by \cite[Exercise II.5.8]{hart}. 
We can extend it to an upper semicontinuous function on
$\mcal{H}\times H$ by formally defining $i_y(C\cap H)$ to be infinity if $C\subset H$.

We shall work with the closed subscheme $W\subset \mcal{H}\times\mcal{H}'$ of 
pairs $(C,X) \in \mcal{H}\times\mcal{H}'$ such that $C$ is a closed
subscheme of $X$ (compare \cite[pp. 265--266]{e-h}, for example).
Let $Z \subset W\times H$ be the closed subscheme of triples
$(C,X,y)$ such that $i_y(C\cap H)\geq e$.  If we choose embeddings 
$\mcal{H}\subset \bfP^N$ and $\mcal{H}'\subset \bfP^{N'}$ into
projective space, then $Z$ is a closed subscheme of 
$\bfP^{N}\times\bfP^{N'}\times H$ defined by finitely many
trihomogeneous polynomials.  
Therefore the fibre $Z_{X,y}$ of the projection from $Z\subset
\mcal{H}\times \mcal{H}'\times H$ to $\mcal{H}'\times H$
is a closed subscheme of $\bfP^{N}$ defined by $O_{d}(1)$ forms of
degree $O_{d}(1)$ in $N+1$ variables.
An application of B\'ezout's theorem (in the form given by 
\cite[Theorem 8.4.6]{fulton}) therefore reveals that
$\#Z_{X,y}=O_d(1)$, provided that $Z_{X,y}$ is zero-dimensional.
This scheme parameterises all closed curves $C\subset X$ of degree $e$
such that $C$ has Hilbert polynomial in $\mcal{P}$ and $i_y(C\cap H)\geq e$.
To complete the proof of Lemma \ref{9cup10}, it is
therefore sufficient to prove that $Z_{X,y}$
is zero-dimensional or empty, whenever $X \cap H$ is integral and $y \in X \cap H$.

We prove this by contradiction.
Let $\mcal{F}_{X,y} \subset X\times
Z_{X,y}$ be the family of closed subschemes of $X$, obtained by
pulling back the universal family $\mcal{F} \subset \bfP^3\times \mcal{H}$
along the projection of $Z\subset
\mcal{H}\times \mcal{H}'\times H$ to $\mcal{H}$.
Let $Y$ be the image of $\mcal{F}_{X,y}$ under the projection onto $X$.  Then $Y$ is a closed
subscheme of $X$, and must have 
dimension at least $2$ if  $\dim Z_{X,y} \geq 1$.
Since $Y\subseteq X$, and $X$ is integral, 
we must therefore conclude that $Y=X$.  
Suppose now that $x \in X \cap H$ is any point
distinct from $y$.  Then there exists a curve
$C\subset X$ of degree $e$, such that $x \in C$ and 
$C$ has intersection multiplicity $e$ with $H$ at $y$.
But then it follows from B\'ezout's theorem that $C\cap H$ contains a
component of dimension $1$, since both $x$ and the $e$-fold point $y$ lie on $C\cap H$.
This plainly contradicts the assumption that $X\cap H$ is integral,
and so completes the proof of Lemma \ref{9cup10}.
\end{proof}

We shall apply Lemma  \ref{9cup10} only in the special cases $e=1$ and $e=2$.
Suppose that $H \in {\bfP^3}^*$ is such that $X\cap H$ is integral,
and let $y \in X\cap H$.  Then in the former case it follows 
from Lemma  \ref{9cup10} that
$$
\#\{\Lambda \in F_1(X): y \in \Lambda\}=O_d(1),
$$
where $F_1(X)$ is the Fano variety $\{\Lambda \in \Grass(1,3): 
\Lambda \subset X\}$ of lines contained in $X$.    In the latter case
it allows us to conclude that 
there are only $O_d(1)$ conics contained in $X$ that are 
tangent to $H$ at $y$.


As above let $X\subset \bfP^3$ be a projective surface.  For 
any non-singular point $x\in X$, we let
$\T_x(X)$ denote the tangent plane to $X$ at $x$.  
The remainder of this section is concerned with 
the set of non-singular points $x\in X$ which have multiplicity at most $2$ on the
intersection $X\cap \T_x(X)$.
We begin by establishing the following result.

\begin{lem}\lab{hess}
Suppose that $X\subset \bfP^3$ is an integral surface which is not
a plane.  Then there exists a non-singular point $x\in X$ such that
$x$ has multiplicity $2$ on the intersection $X\cap \T_x(X)$.
\end{lem}

\begin{proof}
Since $X$ is integral, Bertini's theorem ensures the existence of a 
plane $H\in {\bfP^3}^*$ such that the
intersection $Y=X\cap H$ is integral. 
Now it is well known that given any integral plane curve
of degree exceeding $1$, the Hessian of the curve does not vanish
identically on the curve  (see Fischer \cite[\S 4.5]{fischer}, for
example).  Hence we may find a non-singular point
$x\in Y$ which is not an inflection point.
But then it follows that $x$ is of multiplicity $2$ on the intersection
$X\cap \T_x(X)\cap H$, and thus also on the intersection $X\cap \T_x(X)$.
\end{proof}

The following result allows us to conclude that the set 
of non-singular points $x\in X$ which have multiplicity at most $2$ on the
intersection $X\cap \T_x(X)$, is open in
the Zariski topology.

\begin{lem}\lab{univ}
Let $d$ be a positive integer.  Then there
exists a finite set of universal bihomogeneous polynomials
$$
\Phi_0(a_{\ma{e}}; X_0,X_1,X_2,X_3), \ldots, \Phi_t(a_{\ma{e}};
X_0,X_1,X_2,X_3),
$$
with coefficients in $\Z$, 
such that the following holds for any field $K$.
A point $x=[x_0,x_1,x_2,x_3]$ on the surface defined by the form
\begin{equation}\lab{univ-F}
F=\sum_{\colt{e_0,\ldots,e_3\geq 0}{e_0+e_1+e_2+e_3=d}}
a_{\ma{e}}X_0^{e_0}X_1^{e_1}X_2^{e_2}X_3^{e_3} \in K[X_0,X_1,X_2,X_3],
\end{equation}
is a non-singular point of multiplicity at most $2$ on the section
with the tangent plane at $x$ 
if and only if $\Phi_i(a_{\ma{e}}; x_0,x_1,x_2,x_3)\neq 0$ for some
$0\leq i\leq t$.
\end{lem}

\begin{proof}
Let $X\subset\bfP^3$ be the surface defined by the form
(\ref{univ-F}), and let $x=[\x]\in X$ be a non-singular
point.
Then there exists a linear form $L$ and a quadratic form $Q$,
such that
$$
F(\x+\lambda\y)\equiv 
L(\y)\lambda + Q(\y)\lambda^2  \mod{\lambda^3}.
$$
In fact, if we write $\ma{g}=\nabla F(\x)$ and $\ma{M}=\ma{M}(\x)$ for
the the matrix of second derivatives of $F$
at $\x$, then $L(\y)=\ma{g}.\y$ and $2Q(\y)=\y^{T}\ma{M}\y$.
Moreover it is plain that $\sum_{i=0}^3g_iX_i=0$
is the equation for the tangent plane $\mathbb{T}_x(X)$ at $x$.
It follows that $x=[\x]$ has multiplicity at most $2$ on the tangent
plane section $X\cap\mathbb{T}_x(X)$ if and only if 
$Q(\y)\not=0$ for some vector $\y$ such that $\ma{g}.\y=0$.  
But $\mathbb{T}_x(X)$ is plainly spanned by the
vectors
\begin{align*}
\y_1=(g_1,-g_0,0,0), \quad \y_2=(g_2,0,-g_0,0),\quad 
\y_3=(g_3,0,0,-g_0),\\
\y_4=(0,g_2,-g_1,0),\quad  \y_5=(0,g_3,0,-g_1), \quad
\y_6=(0,0,g_3,-g_2),
\end{align*}
and so  $x\in X$ is a non-singular point of multiplicity at most $2$ on the tangent
plane section if and only if 
$\ma{g}\neq \ma{0}$ and 
$Q(\y_i)\neq 0$
for some $1\leq i\leq 6$.  This therefore establishes the existence of 
the bihomogeneous polynomials 
$$
\Phi_0(a_{\ma{e}}; X_0,X_1,X_2,X_3), \ldots, \Phi_{t}(a_{\ma{e}};
X_0,X_1,X_2,X_3),
$$
that appear in the statement of the lemma.
\end{proof}

Suppose that $X\subset \bfP^3$ is an integral surface defined over a field $K$.
We henceforth let $U$ denote the set of non-singular points on $X$ which have multiplicity at most
$2$ on the tangent plane section at the point. 
On combining Lemma \ref{hess} with Lemma \ref{univ} it follows that 
$U$ is a non-empty open subset of $X$.  Thus 
the complement of $U$ in $X$ consists of finitely many
integral components of dimension at most $1$.  Furthermore, we may
apply B\'ezout's
theorem (in the form given by Fulton \cite[Theorem 8.4.6]{fulton})
to deduce that the sum of the  degrees of these components is bounded in
terms of the maximal degree of the specialisation of the universal
polynomials $\Phi_0, \ldots, \Phi_t$ at $X$.  Thus
the complement of $U$ in $X$ consists of $O_d(1)$ integral components of degree
$O_d(1)$.

\section{Affine surfaces}\lab{aff-surface}

In order to prove Theorem \ref{ash-2} it will be convenient to work
with a projective model for the affine surface $f=0$. Thus let 
$X\subset \bfP^3$ be the surface defined by the form
\begin{equation}\lab{F}
F(X_0,X_1,X_2,X_3)=X_0^\del f(X_1/X_0,X_2/X_0,X_3/X_0),
\end{equation}
of degree $\del$.  
We may assume that $F$ has integer coefficients whose highest common
factor is $1$.  Moreover it is clear that $\ho(f)$ is irreducible if and only if
$X\cap H_\infty$ is integral, where $H_\infty$ denotes the plane $X_0=0$.

As indicated above we shall represent elements of $\bfP^n(\Q)$ by 
$(n+1)$-tuples $\x=(x_0,\ldots,x_n)\in \Z^{n+1}$ such that 
$\hcf(x_0,\ldots,x_n)=1$, for any $n \geq 2.$  The pairs $[\x] =\pm\x$
of such $(n+1)$-tuples corresponds to integral
points  on the scheme $\bfP_\Z^n$.  We shall therefore write 
$[\x] \in \bfP^n(\Z)$ to express that 
$\x\in \Z^{n+1}$ and $\hcf(x_0,\ldots,x_n)=1$.
This choice of notation is motivated by the fact that we shall
constantly deal with integral points on the affine cone over $\Z$,
rather than over $\Q$.  With this in mind we define the counting function 
\begin{equation}\lab{n(B)}
N^\af(\Sigma;B)=\#\{[1,x_1,x_2,x_3] \in \Sigma\cap \bfP^{3}(\Z): 
~H([1,x_1,x_2,x_3]) \leq B\},
\end{equation}
for any locally closed subset $\Sigma\subseteq X$ defined over $\ov{\Q}$.

Let $F \in \Z[X_0,X_1,X_2,X_3]$ be a primitive form  
of degree $d\geq 3$, defining a surface
$X\subset \bfP^3$ such that $X\cap H_\infty$ is integral.
It is now clear that in order to establish Theorem \ref{ash-2} 
it will suffice to establish the estimate
\begin{equation}\lab{lorus}
N^\af(X;B)\ll_{d,\ve} \left\{
\begin{array}{ll}
B^{5/(3\sqrt{3})+1/4+\ve}, & d=3,\\
B^{1+\ve} + B^{3/(2\sqrt{d})+1/3+\ve}, & d\geq 4,
\end{array}
\right.
\end{equation}
for any $\ve>0$.
For any prime $p$ we shall write $X_p$ for the surface defined over $\F_p$ obtained by
reducing the coefficients of $F$ modulo $p$, and we denote 
the set $X_p\cap\bfP^3(\F_p)$ by $X_p(\F_p)$.
Note that, as a scheme, $X_p$ is not necessarily integral.
It will be convenient to define the set 
\begin{equation}\lab{ganga}
S(\Sigma;B)=\{[1,x_1,x_2,x_3] \in \Sigma\cap \bfP^{3}(\Z): 
~H([1,x_1,x_2,x_3]) \leq B\},
\end{equation}
for any locally closed subset $\Sigma\subseteq X$ defined over $\ov{\Q}$.
In particular we have $N^\af(\Sigma;B)=\#S(\Sigma;B)$, by (\ref{n(B)}).
Now let $\pi=[1,\pi_1,\pi_2,\pi_3]\in X_p(\F_p)$, where 
$\pi_1,\pi_2,\pi_3$ are always assumed to be in $\F_p$.
We also define the set 
$$
S_p(\Sigma;B,\pi)=\left\{[1,x_1,x_2,x_3] \in \Sigma\cap \bfP^{3}(\Z): 
\begin{array}{l}
H([1,x_1,x_2,x_3]) \leq B,\\ 
x_i\equiv \pi_i \mod{p}, ~(1\leq i \leq 3)
\end{array}
\right\},
$$
for any locally closed subset $\Sigma\subseteq X$ defined over $\ov{\Q}$.
It is clear that 
$S_p(\Sigma;B,\pi)$ may be empty if there are no points on $\Sigma$
which specialise to $\pi$ on $X_p$.

Recall the definition of the non-empty open subset $U\subset X$,
introduced at the close of \S \ref{geometry}.
For any prime $p$ we shall define $U_p$ to be
the corresponding open set of non-singular
points on $X_p$ which have multiplicity at most $2$ on the tangent plane
section at the point.
Our main tool in the proof of Theorem \ref{ash-2} will be the following
adaption of a result due to the second author \cite[Theorem 14]{annal}.

\begin{lem}\lab{14}
Let $\ve>0$ and let 
$X\subset \bfP^3$ be an integral surface, defined by a
primitive form  $F \in \Z[X_0,X_1,X_2,X_3]$
of degree $d\geq 3$ such that 
\begin{equation}\lab{height}
\log H(F)=O_{d}(\log B).
\end{equation}
Then there exists 
a set $\Pi$ of $O_{d,\ve}(1)$ primes $p$, with
\begin{equation}\lab{prime}
B^{1/{\sqrt{d}}+\ve}\ll_{d,\ve} p\ll_{d,\ve} B^{1/{\sqrt{d}}+\ve},
\end{equation}
such that the following holds.
For each $\pi=[1,\pi_1,\pi_2,\pi_3]\in U_p(\F_p)$, there exists a form
$G_\pi \in \Z[X_0,X_1,X_2,X_3]$  
of degree $O_{d,\ve}(1)$ which is not divisible by $F$, such that
$$
S(U;B)=
\bigcup_{p\in \Pi}\bigcup_{\pi\in U_p(\F_p)}
\{[1,x_1,x_2,x_3] \in
S_p(U;B,\pi): G_\pi(1,x_1,x_2,x_3)=0\}.
$$
\end{lem}

\begin{proof}
The proof of Lemma \ref{14} is based upon an application of
\cite[Theorem 14]{annal} in the case $n=4$ and $\mathbf{B}=(1,B,B,B)$.
Let $p$ be any prime, with $p\geq CB^{1/{\sqrt{d}}+\ve}$ for an
appropriate constant $C$ depending at most upon $d$ and $\ve$.  Furthermore, let 
$\pi=[1,\pi_1,\pi_2,\pi_3]\in X_p(\F_p)$ be any non-singular point on the
reduction $X_p$ of $X$ modulo $p$.
Then a  rudimentary examination of the proof 
shows that there exists a form 
$G_\pi \in \Z[X_0,X_1,X_2,X_3]$  
of degree $O_{d,\ve}(1)$ which is not divisible by $F$, such that
$G_\pi$ vanishes at all points of the set $S_p(X;B,\pi)$.

Our task is to show that there
exists a set $\Pi$ of $O_{d,\ve}(1)$ primes $p$ which all satisfy (\ref{prime}),
such that for each $x\in S(U;B)$ there
exists $p \in \Pi$ for which the reduction $\pi$ of the
point $x$ modulo $p$ lies in the open subset $U_p$.
This will clearly suffice to establish the lemma.
For $0\leq i \leq t$ let 
$$
\phi_i(X_0,X_1,X_2,X_3)=\Phi_i(a_{\ma{e}};
X_0,X_1,X_2,X_3) \in \Z[X_0,X_1,X_2,X_3]
$$ 
be the forms obtained by
specialising the universal polynomials in Lemma \ref{univ} at the
primitive form $F$.  It is clear that $t=O_d(1)$ and that 
the degree of each $\phi_i$ is $O_d(1)$.  Moreover, it also follows that
$\log H(\phi_i)=O_d(\log B)$ for $0\leq i\leq t$, by (\ref{height}).
Now let $x\in S(U;B)$.  Then by the remark after Lemma \ref{univ} we
may assume that there exists $0\leq i\leq t$ such that
$\phi_i(1,x_1,x_2,x_3)\neq 0$.  In particular it is clear that 
$$
\log |\phi_i(1,x_1,x_2,x_3)|\ll_d \log B.
$$ 
Hence the number of primes $p\geq
CB^{1/{\sqrt{d}}+\ve}$ such that $p\mid \phi_i(1,x_1,x_2,x_3)$ is at
most $c(d,\ve)$.
By Bertrand's postulate we may select a set $\Pi$ of primes satisfying
(\ref{prime}) such that $\#\Pi > c(d,\ve)$.
Thus $\Pi$ must contain a prime $p$ for which $p\nmid \phi_i(1,x_1,x_2,x_3)$. 
We complete the proof
of Lemma \ref{14} by applying Lemma
\ref{univ} with the choice of field $K=\F_p$.
\end{proof}

It is clear from Lemma \ref{14} that we shall need to consider the 
contribution to $N^\af(X;B)$ arising from points lying on
a finite set of curves that are contained in the
surface $X$.  
Let $I_{D}$ be any finite set of integral curves $C$ contained in $X$, 
each of degree $D$, and define
$$
\Sigma_D=\bigcup_{C\in I_D}C.
$$
The following result is a crucial ingredient in the
proof of Theorem \ref{ash-2}, and will be established in the following
subsection.  In fact it can be thought of as providing independent evidence for
\conj \ref{con-ash}.

\begin{pro}\lab{burn}
Let $\ve>0$ and 
suppose that $X\subset \bfP^3$ is a surface of degree $d\geq 3$ such that
$X\cap H_\infty$ is integral.
Then we have
$$
N^\af(\Sigma_D;B) \ll_{d,\ve} B^{\ve}\max\{B^{2/d},B^{1/D},\#I_{D}\}
$$
for $D=1$ or $2$.
\end{pro}

\subsection{Proof of Proposition \ref{burn}}\lab{pf-main'}

Let $X\subset \bfP^3$ be a surface of degree $d\geq 3$ such that
$X\cap H_\infty$ is integral, and let $I_D$ be a finite set of curves
$C$ contained in $X$, each of degree $D$.  
We begin with the case $D=1$.  
Any line $L\in I_1$  which contains at most one 
point  $[1,x_1,x_2,x_3] \in \bfP^{3}(\Z)$ such that
$H([1,x_1,x_2,x_3]) \leq B$, 
clearly contributes $O(1)$ to $N^\af(\Sigma_1;B)$.  Such lines are
therefore satisfactory from the point of view Proposition \ref{burn}.  
Suppose now that $L \in I_1$ contains more than one
rational point $[1,x_1,x_2,x_3]$ of height at most $B$.  We choose two
such points $[1,\ma{t}]$ and $[1,\ma{t}+\ma{s}]$, 
for $\ma{t}, \ma{s}\in \Z^3$ such that $|\ma{s}|$ is minimal.
In the notation of (\ref{ganga}), it is then clear that every member of the set $S(L;B)$ is represented
by an integer vector of the form 
$(1,\bt+n\ma{s})$, for some $n\in\Z$. 
In particular it follows from Lemma \ref{4count} that 
\begin{equation}\lab{coffee}
N^\af(L;B)=\#S(L;B)=O(1+B/|\ma{s}|).
\end{equation}
Now the vector $\ma{s}$ produces a rational
point $[0,\ma{s}]$ on the integral curve $X\cap
H_{\infty}$, and by Lemma \ref{9cup10} there are $O_d(1)$ lines through
such a point.  
Thus the contribution from all the lines under 
consideration passing through $[0,\ma{s}]$ is $O_d(1+B/|\ma{s}|)$.  It
remains to sum over all rational points $[0,\ma{s}]\in X\cap H_{\infty}$
for which $|\ma{s}|\ll B$.  We do this by considering separate dyadic
ranges $S<|\ma{s}|\leq 2S$.  According to Lemma \ref{bub} there are at
most 
$$
N_{X\cap H_\infty}(2S)\ll_{d,\ve}S^{2/d+\ve}
$$ 
such points. Thus the overall
contribution is
$O_{d,\ve}(S^{2/d+\ve}+BS^{2/d-1+\ve})$ from this range. When we sum this up with
$S\ll B$ running over powers of $2$ we obtain a total
$O_{d,\ve}(B^{1+\ve})$, as required.  This completes the argument for
the case $D=1$.

Turning to the case $D=2$, we first record a basic estimate for the
number of rational points of bounded height on certain conics, in
which one of the coordinates is fixed.
Let $C\subset \bfP^2$ be a plane conic defined by a 
non-singular primitive quadratic form $q\in \Z[X_0,X_1,X_2]$.
Suppose that the coefficients of $q$ are bounded
in modulus by $H(q)$, and that the binary form 
$q(0,X_{1},X_{2})$ is also
non-singular.  
Then the second author has shown that
$$
\#\{[k,x_1,x_2] \in C\cap \bfP^{2}(\Z): 
~H([k,x_1,x_2]) \leq B\}\ll_\ve H(q)^{\ve}B^{\ve}
$$
for any fixed integer $k$ and any choice of $\ve>0$ \cite[Theorem 3]{cubic}.
We can remove the dependence of this estimate upon $H(q)$
by applying Lemma \ref{theta}.  Indeed 
we find that either
there are $O(1)$  points, or else $H(q)\ll B^{12}$.  In either case the
total number of solutions is $O_{\ve}(B^{13\ve})$.  We may therefore conclude as
follows.

\begin{lem}\label{affquad2}
Let $\ve>0$ and suppose that
$C\subset \bfP^2$ is a plane conic defined by a 
non-singular quadratic form $q\in \Z[X_0,X_1,X_2]$ such that 
the binary form 
$q(0,X_{1},X_{2})$ is also
non-singular.  
Then for any integer $k$ we have
$$
\#\{[k,x_1,x_2] \in C\cap \bfP^{2}(\Z): 
~H([k,x_1,x_2]) \leq B\}\ll_\ve B^{\ve}.
$$
\end{lem}

We are now ready to handle those conics $C\in I_2$  
which are not tangent to the plane $H_{\infty}\in {\bfP^3}^*$.
Such a conic is specified by a linear equation
$a_0X_0=a_1X_1+a_2X_2+a_3X_3$, in which not all the coefficients
$a_1,a_2,a_3$ vanish, together with
a homogeneous quadratic equation $Q(X_0,X_1,X_2,X_3)=0$.  
If $a_3\not=0$,
say, we may use the
linear equation to eliminate $X_3$, so as to produce an equation
$q(X_0,X_1,X_2)=0$.  
Since $C$ is assumed to be integral it 
follows that $q$ must be non-singular.  Similarly, since $C$ is not
tangent to $H_{\infty}$ the binary form $q(0,X_1,X_2)$ will also be
non-singular. We may therefore apply Lemma \ref{affquad2} to deduce that
$$
\#\{[1,x_1,x_2,x_3] \in C\cap \bfP^{3}(\Z): ~H([1,x_1,x_2,x_3]) \leq B\}\ll_\ve B^{\ve}.
$$
Thus conics which are not tangent
to $H_{\infty}$ make a satisfactory contribution to $N^\af(\Sigma_2;B)$.

It remains to consider conics $C$ which are tangent to $H_{\infty}$.
We plan to show that the affine integral points on such a conic 
can be parameterised by the values at integer points of a small
number of quadratic
polynomials.  To estimate how many such
polynomials are necessary we must first bound the coefficients of the
defining equations for $C$.  As above $C$ can be given by a linear equation
$a_0X_0=a_1X_1+a_2X_2+a_3X_3$, in which $a_3$, say, is non-zero,
and a homogeneous quadratic equation in which we can eliminate $X_3$
to produce the defining equation
$q(X_0,X_1,X_2)=0$ for $C$.   
An application of Lemma \ref{theta} shows that either $C$ contains $O(1)$
rational points of height at most $B$, or else $H(q)\ll B^{12}$, as we
henceforth assume.  We proceed by showing that the coefficients of the
linear form can be taken to satisfy
\begin{equation}\lab{ai}
|a_0|\ll B^3, \quad |a_1|,|a_2|,|a_3|\ll B^2.
\end{equation}
To do so we first consider the case 
in which the points 
$[1,x_1,x_2,x_3]\in\bfP^3(\Z)$ of height at most $B$, which lie on the plane
defined by the linear equation, are restricted to a line contained in
the plane.   In this case the relevant points on $C$ are
also restricted to a line, so that there are at most two rational points
to be counted.  In the alternative case we can find  
rational points $[1,x_1^{(i)},x_2^{(i)},x_3^{(i)}]$ contained 
in the plane, for $1\leq
i \leq 3$, which lie in general 
position and have height at most $B$.
It follows that one can write a suitable scalar multiple
of the vector $(a_0,a_1,a_2,a_3)$ in terms of determinants in 
the integer vectors 
$(1,x_1^{(i)},x_2^{(i)},x_3^{(i)})$.
Taking ${\rm h.c.f.}(a_0,a_1,a_2,a_3)=1$ this establishes (\ref{ai}).

Next we obtain a preliminary polynomial parameterisation of $C$.  Since $C$
is tangent to $H_{\infty}$ the binary form $q(0,X_1,X_2)$ must have 
rank $1$.  This allows us to write 
$$
q(X_0,X_1,X_2)=a(\alpha X_1+\beta X_2)^2+bX_0X_1+cX_0X_2+dX_0^2,
$$
for integers $a,b,c,d,\alpha,\beta\ll B^{12}$ with $\al, \be$ 
coprime.  Choose integers
$\gamma,\delta\ll B^{12}$ with $\al\del-\be\gamma=1$ and set
\begin{equation}\label{subst}
X_0=Y_0, \quad \alpha X_1+\beta X_2=Y_1,\quad \gamma X_1+\del X_2=Y_2,
\end{equation}
so that
$$
q(X_0,X_1,X_2)=q'(Y_0,Y_1,Y_2)=aY_1^2+eY_0Y_1+fY_0Y_2+dY_0^2
$$
for integers $e,f\ll B^{24}$.  Clearly the zero locus of $q$ coincides
precisely with the zero locus of $q'$.  
If $f$ were to vanish the conic $C$ would degenerate to a 
pair of lines, which is
a case we have excluded.  Hence $f\not=0$, so that if $X_0=Y_0=1$
then
$$
Y_2=-f^{-1}(aY_1^2+eY_1+d)
$$
and hence
$$
X_1=\del Y_1-\be Y_2=q_1(Y_1),\;\;\; X_2=\al Y_2-\gamma Y_1=
q_2(Y_1)$$
and
$$ 
X_3=a_3^{-1}(a_0-a_1X_1-a_2X_2)=q_3(Y_1)
$$
where $q_1,q_2,q_3$ are quadratic polynomials with rational
coefficients whose common denominator is $O(B^\kappa)$
for some absolute constant $\kappa$. 

Our problem is now to estimate the number of integers $Y$ for which
$$
(q_1(Y),q_2(Y),q_3(Y))\in\Z^3, \quad
\max\{|q_1(Y)|, |q_2(Y)|,|q_3(Y)|\} \leq B.
$$
If there
is no such integer we are clearly done.  Otherwise let $Y^{(*)}$ be any
such integer and substitute $Y=Y^{(*)}+Z$ to produce polynomials
$Q_i(Z)$ each with an integer constant term.  We then write $D$ for
the lowest common denominator of the coefficients of the polynomials
$Q_i$ so that
$$
Q_i(Z)=A_i+D^{-1}(B_iZ+C_iZ^2),\quad A_i, B_i, C_i\in\Z,
$$
for $i=1,2,3$.
In view of our estimate for the common denominator of the polynomials
$q_i$ we have 
\begin{equation}\label{Dkappa}
D\ll B^{\kappa}.
\end{equation}
The polynomials $Q_i$ provide our ``preliminary parameterisation''.

We proceed to examine the set of integers $Z$ for which
$Q_1(Z),Q_2(Z),Q_3(Z)$ are all integral.  We classify such $Z$
according to the value of ${\rm h.c.f.}(Z,D)=\la$, say, and we write
$D=\la\mu$. The number of possible
classes is $d(D)\ll_{\ve}B^{\ve}$ in view of (\ref{Dkappa}).  It
is precisely for this estimate that it is necessary to control the
size of the coefficients in our preliminary parameterisation.  We now
claim that if the class corresponding to $\la$ is non-empty 
then there exists $Z_{\la}\in\Z$ and $D_{\la}|D$ such that
\begin{align*}
\{Z\in\Z:\,(Q_1(Z),Q_2(Z),Q_3(Z))\in\Z^3,\; &{\rm h.c.f.}(Z,D)=\la\}\\
&=\{Z\in\Z:\,Z\equiv Z_{\la}\mod{D_{\la}}\}.
\end{align*}
Once this claim is established we will write 
$$
R_{i,\la}(t)=Q_i(Z_{\la}+D_{\la}t).
$$
Then the quadratic polynomials $R_{i,\la}(t)$ are integer valued, so
that $2R_{i,\la}(T)\in\Z[T]$.  Moreover any integer point on our
affine conic will be of the form
$[1,R_{1,\la}(t),R_{2,\la}(t),R_{3,\la}(t)]$ with $t\in\Z$ for some
value of $\la$.  This will
produce the required set of integer parameterisations, one for each
value of $\la$.

We must now verify the above claim.  Write $Z=\la W$.  Then the condition
$Q_i(Z)\in\Z$ is equivalent to $\mu|B_iW+C_i\la W^2$. However the relation
${\rm h.c.f.}(Z,D)=\la$ implies that $\mu$ and $W$ are coprime, so
that in fact $Q_i(Z)\in\Z$ if and only if $\mu|B_i+C_i\la W$.  We now
write $D_i={\rm h.c.f.}(\mu,C_i\la)$ and $\mu=D_i\mu_i$.  Then for
each $i$ the set of solutions $W$ of the congruence $C_i\la W\equiv
-B_i\mod{\mu}$ is either empty (if $D_i\nmid B_i$), or consists of a single
residue class $Z_{i,\la}\mod{\mu_i}$.  To complete the proof of the claim it
suffices to observe that the intersection of the three residue classes
$Z\equiv Z_{i,\la}\mod{\mu_i}$ is either empty, or is a single residue
class modulo ${\rm l.c.m}(\mu_1,\mu_2,\mu_3)$.

We can at last complete our treatment of the case
$D=2$.  Consider those rational points $[1,x_1,x_2,x_3]\in C\cap \bfP^3(\Z)$
which are given by 
$$ 
x_i=A_i+B_it+C_it^2, \quad (i=1,2,3),
$$
for $t\in \Z$ and fixed $A_i,B_i,C_i\in \frac{1}{2}\Z$.
Since these polynomials parameterise $C$ it follows that the point
$[0,C_1,C_2,C_2]$ lies on $C\cap H_{\infty}$.  Moreover, Lemma 
\ref{4count} shows that there are 
$$
\ll 1+\lb\frac{B}{\max\{|C_1|,|C_2|,|C_3|\}}\rb^{1/2}
$$ 
corresponding points belonging to the set $S(C;B)$, in the 
notation of (\ref{ganga}).
Now let $p=[0,a,b,c]\in X\cap H_{\infty}$, with ${\rm
  h.c.f.}(a,b,c)=1$, 
and suppose that 
$S<\max\{|a|,|b|,|c|\}\leq 2S$.  Then if $C\subset X$ is
tangent to $H_{\infty}$ at $p=[0,a,b,c]$ it follows that the
contribution from $C$ to $N^\af(\Sigma_2;B)$ is
$O_{\ve}(B^{\ve}(1+B^{1/2}S^{-1/2}))$, after
allowing for the various parameterisations of $C$.  According to Lemma
\ref{9cup10} we get the overall bound $O_{d,\ve}(B^{\ve}(1+B^{1/2}S^{-1/2}))$
when we include all other conics in $X$ which are 
tangent to $H_{\infty}$ at $p$.  It remains to sum 
over all rational points $p$ on $X\cap H_{\infty}$
for which $\max\{|a|,|b|,|c|\}\ll B$.  We do this by the same dyadic
range decomposition as was used for the case $D=1$.  Each range of
length $S$ will contribute
$$
\ll_{d,\ve}B^{\ve}(S^{2/d+\ve}+B^{1/2}S^{2/d-1/2+\ve}),
$$ 
giving a
total contribution $O_{d,\ve}(B^{2\ve}(B^{2/d}+B^{1/2}))$ to
$N^\af(\Sigma_2;B)$. This is plainly satisfactory for Proposition
\ref{burn}.

\subsection{Completion of the proof of Theorem \ref{ash-2}}\lab{pf-main}

Let $F \in \Z[X_0,X_1,X_2,X_3]$ be a primitive form  
of degree $d\geq 3$, defining an integral surface
$X\subset \bfP^3$.
Let $H_\infty\in {\bfP^3}^*$ denote the plane $X_0=0$.
In this section we shall establish the estimate (\ref{lorus}), and so
complete the proof of Theorem~\ref{ash-2}.  We therefore 
assume henceforth that 
the intersection $X\cap H_\infty$ is integral.
Suppose first that the coefficients of $F$ are large
compared with $B$.  Then an application of Lemma \ref{theta} reveals
that $N^\af(X;B)\leq N^\af(Y;B)$, for some curve $Y\subset \bfP^3$ obtained by
intersecting $X$ with a distinct surface $F'=0$.
The contribution from irreducible components of $Y$ of degree 
$2$ or more is $O_{d,\ve}(B^{1+\ve})$, by Lemma~\ref{bub}, while for 
components of degree $1$ the bound (\ref{coffee}) again gives a  
satisfactory result. We henceforth assume that  (\ref{height}) holds for $F$.

Recall the definition of the non-empty open subset 
$U\subset X$, introduced at the close of \S \ref{geometry}. 
Then by the comments made directly after Lemma \ref{univ} we may
argue as above to deduce that
$N^\af(X\setminus U;B)\ll_{d,\ve} B^{1+\ve}$.  Hence it is 
sufficient  to establish the estimate 
\begin{equation}\lab{goal}
N^\af(U;B)\ll_{d,\ve} \left\{
\begin{array}{ll}
B^{5/(3\sqrt{3})+1/4+\ve}, & d=3,\\
B^{1+\ve} + B^{3/(2\sqrt{d})+1/3+\ve}, & d\geq 4.
\end{array}
\right.
\end{equation}
We are now in a position to apply Lemma \ref{14}.
Thus there exists
a set $\Pi$ of $O_{d,\ve}(1)$ primes $p$, with
(\ref{prime}) holding, such that 
$$
N^\af(U;B)\leq
\sum_{p\in \Pi}
\sum_{\pi\in U_p(\F_p)}
\hspace{-0.2cm}
\#\{[1,x_1,x_2,x_3] \in
S_p(U;B,\pi): G_\pi(1,x_1,x_2,x_3)=0\}.
$$
Here the $G_\pi \in \Z[X_0,X_1,X_2,X_3]$  are a finite set of forms
indexed by points $\pi=[1,\pi_1,\pi_2,\pi_3]\in U_p(\F_p)$. 
For each $\pi\in U_p(\F_p)$, the form 
$G_\pi$ has degree $O_{d,\ve}(1)$ and is not divisible by $F$.

Let $Y=Y_\pi \subset \bfP^3$ be an integral component of the curve
$F=G_\pi=0$, and let 
$$
N_p^\af(Y;B,\pi)=\#S_p(Y;B,\pi).
$$
If $Y_1\cup\cdots\cup Y_s$ denotes the irreducible decomposition
of $F=G_\pi=0$,
with its reduced scheme structure, then we clearly have
$s=O_{d,\ve}(1)$ and
$$
\#\{[1,x_1,x_2,x_3] \in
S_p(U;B,\pi): G_\pi(1,x_1,x_2,x_3)=0\}
\leq \sum_{1\leq j\leq s} N_p^\af(Y_j;B,\pi).
$$
As we vary over primes $p\in \Pi$ and points 
$\pi\in U_p(\F_p)$, let $I$ denote the set of integral components 
of the curves $F=G_\pi=0$ which have degree at most $2$.  
Since $\#U_p(\F_p)\leq
\#X_p(\F_p) = O_d(p^2)$, it is clear from (\ref{prime}) that
$$
\mbox{$\#I\leq \#\Pi.\max_p \#U_p(\F_p)= O_{d,\ve}(B^{2/\sqrt{d}+2\ve})$}. 
$$
Hence Proposition \ref{burn} implies that the overall contribution
to $N^\af(U;B)$ from the set $I$ is 
$$
\ll_{d,\ve}B^{1+\ve}+ B^{2/\sqrt{d}+3\ve}.
$$
This is satisfactory for (\ref{goal}) provided that $d\geq 3$.

We henceforth fix a choice of prime $p\in \Pi$, and a point 
$\pi\in U_p(\F_p)$.  
Let $G_\pi$ be the corresponding form produced by
Lemma \ref{14}, and let 
$Y=Y_\pi \subset \bfP^3$ be any integral component of the curve
$F=G_\pi=0$, of degree $e$. It henceforth suffices to assume that
$Y$ has dimension $1$ and degree $e\geq 3$, since we have already
taken care of the remaining components.  In particular we recall that
$e=O_{d,\ve}(1)$.   
The main work in this section will be taken up with
establishing the following estimate. 

\begin{pro}\lab{pro}
Assume that $Y$ has degree $e\geq 3$.  Then we have 
$$
N_p^\af(Y;B,\pi)\ll_{d,\ve} B^{1/e-1/((e-1)\sqrt{d})}.
$$
\end{pro}

When $d\geq 4$ the estimate in Proposition \ref{pro} takes its largest
value at $e=3$, whereas when $d=3$ it is maximal at $e=4$.  
We therefore have
$$
N_p^\af(Y;B,\pi)\ll_{d,\ve} \left\{
\begin{array}{ll}
B^{1/4-1/(3\sqrt{3})}, & d=3,\\
B^{1/3-1/(2\sqrt{d})}, & d\geq 4,
\end{array}
\right.
$$
whenever $Y$ has degree $e\geq 3$.
On inserting this bound 
into our previous argument and re-defining $\ve$, we therefore arrive at the inequality
(\ref{goal}). Hence it will suffice to prove Proposition \ref{pro} in
order to complete the proof of Theorem \ref{ash-2}.

We shall prove Proposition \ref{pro} by using a generalisation
of the determinant method, as developed by the second author  \cite[Theorem 14]{annal}.
By a result of Gruson, Lazarsfeld and Peskine \cite{mum}, the ideal of $Y$ is generated by forms $F_1,\ldots,F_r$
whose maximal degree is bounded in terms of $e$ alone.
We may also assume  that these forms are linearly independent, and that
the number of forms is bounded in terms of $e=O_{d,\ve}(1)$.
By applying Broberg's generalisation \cite[Lemma 5]{broberg} of Lemma \ref{theta}, we may
assume that either $\log H(F_i) =O_{d,\ve}(\log B)$ for $1\leq i\leq r$,
or else 
$$
N_p^\af(Y;B,\pi)=O_{d,\ve}(1),
$$ 
which is satisfactory for the statement of Proposition \ref{pro}.
We henceforth assume that  the forms defining $Y$ 
all have logarithmic height $O_{d,\ve}(\log B)$.

Let $q$ be a prime, and let $J=(F_1, \ldots, F_r)\subset \Q[X_0,X_1,X_2,X_3]$ be the
homogeneous prime ideal defining $Y$.
We write $Y_q=\rom{Proj}(\F_q[X_0,X_1,X_2,X_3]/J_q)$, where $J_q$
denotes the image of $J\cap \Z_{(q)}[X_0,X_1,X_2,X_3]$ in 
$\F_q[X_0,X_1,X_2,X_3]$.  Here, as throughout this paper, $\Z_{(q)}$
denotes the localisation of $\Z$ at the prime $q$.
Let $V\subset Y$ be the non-empty open subset 
of non-singular points on the curve $Y$. 
We shall write $V_q$ for the set of non-singular points on 
$Y_q$. Finally, for any choice of $\om=[1,\om_1,\om_2,\om_3]\in V_q(\F_q)$ we 
set
$$
S_{p,q}(V;B,\pi,\om)=\left\{x \in 
S_{p}(V;B,\pi): x_i\equiv \om_i \mod{q}, ~(1\leq i \leq 3)
\right\}.
$$
We recall here that $x=[1,x_1,x_2,x_3]$ for $x_i\in \Z$ such
that $x_i\equiv \pi_i \mod{p}$ for $1\leq i\leq 3$,
whenever $x \in  S_{p}(V;B,\pi)$.
We proceed by establishing the following result.

\begin{lem}\lab{q}
Let $a>0$.  Then there exists a set $\Omega$ of $O_{d,\ve}(1)$ primes
$q$, with $p\not\in \Omega$ and 
$$
B^{a} \ll_{d,\ve} q \ll_{d,\ve} B^{a},
$$
such that 
$$
S_p(V;B,\pi)=
\bigcup_{q\in \Omega}\bigcup_{\colt{\om\in
    V_q(\F_q)}{\om=[1,\om_1,\om_2,\om_3]}} S_{p,q}(V;B,\pi,\om).
$$
\end{lem}

\begin{proof}
The proof of Lemma \ref{q} is straightforward.  Recall that the ideal
$J$ of $Y$ is generated by $O_{d,\ve}(1)$ forms of maximal degree $O_d(1)$,
all of which have logarithmic height $O_{d,\ve}(\log B)$.
Using the Jacobian
criterion for a point of $Y$ to be contained in the singular locus, we
may find $O_{d,\ve}(1)$ non-zero forms  
$$
\Psi_1, \ldots, \Psi_l\in \Z[X_0,X_1,X_2,X_3]
$$  
of maximal degree $O_{d,\ve}(1)$, such that 
$\Psi_i(1,x_1,x_2,x_3)\neq 0$ for some $1\leq i \leq l$ whenever
$x=[1,x_1,x_2,x_3] \in S_p(V;B,\pi)$.
On noting that 
$$\log |\Psi_i(1,x_1,x_2,x_3)|\ll_{d,\ve} \log B$$ 
and
recalling from (\ref{prime}) that $\log p\ll_{d,\ve} \log B$, one 
easily deduces the statement of Lemma \ref{q} along precisely the same
lines as those used to prove Lemma \ref{14}.
\end{proof}

It follows from the proof of Lemma \ref{q} that $Y\setminus V$ is a
variety of dimension $0$ and degree $O_{d,\ve}(1)$, whence
$N_p^\af(Y\setminus V;B,\pi)=O_{d,\ve}(1)$.  In order to prove Proposition
\ref{pro} it therefore suffices to establish the estimate
$$
N_p^\af(V;B,\pi)\ll_{d,\ve} B^{1/e-1/((e-1)\sqrt{d})}.
$$
To do so we shall apply Lemma \ref{q} with $a=1/e-1/((e-1)\sqrt{d})$,
so that in particular 
\begin{equation}\lab{q-prime}
B^{1/e-1/((e-1)\sqrt{d})} \ll_{d,\ve} q \ll_{d,\ve} B^{1/e-1/((e-1)\sqrt{d})},
\end{equation}
for every $q \in \Omega$.  On noting that $\#V_q(\F_q)\leq \#Y_q(\F_q)=O_d(q)$, it plainly 
suffices to show that 
\begin{equation}\lab{end-game}
S_{p,q}(V;B,\pi,\om)\ll_{d,\ve} 1,
\end{equation}
for each $q\in \Omega$ and each corresponding point $\omega\in V_q(\F_q)$.

We henceforth fix a choice of $q\in \Omega$,  and a point 
$\om=[1,\om_1,\om_2,\om_3]\in V_q(\F_q)$.
We shall establish the existence of an auxiliary form, not contained
in the prime ideal $J$ defining $Y$, which vanishes at all points 
$S_{p,q}(V;B,\pi,\om)$.  An application of B\'ezout's theorem will
then yield the estimate (\ref{end-game}), provided we can show that
the degree of this auxiliary form is bounded in terms of $d$ and
$\ve$ alone.  

For a fixed positive integer $k$, we shall need to find a 
set of monomials $M_1,\ldots,M_k\in\Z[X_0,X_1,X_2,X_3]$ of degree
$D$ depending only on $e$ and $k$,  
such that no non-trivial linear combination of them belongs to the
prime ideal $J$, and such that the sum
$$
\deg M_1(1,T_1,T_2,T_3) +\cdots +\deg M_k(1,T_1,T_2,T_3)
$$
is as small as possible.  We shall henceforth write $m_i(T_1,T_2,T_3)$ to
denote the monomial $M_i(1,T_1,T_2,T_3)$, for $1\leq i \leq k$.
The following result shows that such a set always exists provided that
$k$ is taken to be sufficiently large.

\begin{lem}\lab{upper}
Let $k\gg_e 1$ be a positive integer.   
Then there exists a set  $\mcal{M}(J,k)$ 
of monomials $M_1,\ldots,M_k\in \Z[X_0,X_1,X_2,X_3]$ of degree
$D=O_{d,k}(1)$,  such that no non-trivial linear combination of them 
belongs to the
prime ideal $J$, and such that 
\begin{equation}\lab{matharan}
\deg m_1 +\cdots +\deg m_k \leq
\frac{k^2}{2e} +O_e(k). 
\end{equation}
\end{lem}

\begin{proof}
Let $H_\infty\in {\bfP^3}^*$ be the plane $X_0=0$, and let $Z$ be the
scheme-theoretic intersection of $Y$ and $H_\infty$.  Then there exists an
integer $\del_0$ depending only on $e$, such that the Hilbert function
$h_Z(\del)$ of $Z\subset \bfP^3$ satisfies $h_Z(\del)=e$ for $\del\geq
\del_0$. Whenever $\del\geq
\del_0$ we may therefore find $e$ monomials in $\Z[X_1,X_2,X_3]$ of
degree $\del$, such that no non-trivial linear combination of them
belongs to the ideal in $\Q[X_0,X_1,X_2,X_3]$ which is generated 
by $J$ and $X_0$.
Denote these monomials by $m_{1,\del}, m_{2,\del},\ldots,m_{e,\del}$.
We then examine the monomials
$$
M_{i,\del}=X_0^{D-\del}m_{i,\del}(X_1,X_2,X_3),\;\;\; (1\le i\le
e,\;\; \del_0\le\del\le D).
$$
Suppose that some linear
combination of them belonged to $J$.  We would then have a relation of
the form
$$
\sum_{\del_0\le\del\le D}X_0^{D-\del}R_{\del}(X_1,X_2,X_3)\in J,
$$
where $R_\del$ is a linear combination of the monomials
$m_{1,\del}, m_{2,\del},\ldots,m_{e,\del}$.  By our construction, if
$R_{\del}\in \langle J,X_0\rangle$ then the coefficients of
$m_{1,\del},\ldots, m_{e,\del}$ must all vanish, so that $R_{\del}=0$.  Using this observation for $R_D,
R_{D-1},R_{D-2},\ldots$ we find that no non-trivial linear combination
of the $M_{i,\del}$ can belong to $J$.

Writing $n_D$ for the number of monomials under consideration we have 
$n_D=eD+O_e(1)$.  Moreover the sum
of the degrees of the corresponding $m_{i,\del}$ is $eD^2/2+O_e(D)$.
The result then follows on choosing $D$ to be  the smallest integer for which
$n_D\ge k$, and ignoring any surplus monomials.
\end{proof}

Let $k\gg_e 1 $ be a positive integer and suppose that
$\#S_{p,q}(V;B,\pi,\om)\geq k$.  We select a list of points 
$$
x^{(1)}, \ldots, x^{(k)}\in S_{p,q}(V;B,\pi,\om).
$$ 
Then, for each $1\leq j\leq k$, we have 
$x^{(j)}=[1,x_1^{(j)},x_2^{(j)},x_3^{(j)}]$ for $x_i^{(j)}\in \Z$ such
that 
\begin{equation}\lab{reduc}
x_i^{(j)}\equiv \pi_i \mod{p}, \quad x_i^{(j)}\equiv \om_i \mod{q},
\quad (i=1,2,3).
\end{equation}
We define $\D$ to be the determinant of the
$k\times k$ matrix
$$
\big(M_i(1,x_1^{(j)},x_2^{(j)},x_3^{(j)})\big)_{1\leq i,j \leq k},
$$
where $M_1,\ldots,M_k\in \mcal{M}(J,k)$, so that in particular we have
$$
\D=\det\big(m_i(x_1^{(j)},x_2^{(j)},x_3^{(j)})\big)_{1\leq i,j \leq k}.
$$
Thus it follows from an application of Lemma \ref{upper} 
that
\begin{equation}
\log |\D| \leq  k\log k +\frac{k^2}{2e}\log B +O_e(k\log B).
\lab{D1}
\end{equation}

We shall also need to show that the determinant $\D$ is divisible by large powers of $p$
and $q$.  Once combined with (\ref{D1}), this information  will then 
be used to establish that $\D=0$ provided that $k$ is chosen to be
sufficiently large in terms of $d$ and $\ve$.  
In particular, once such a value of $k$ is fixed, it will follow that
the monomials in $\mcal{M}(J,k)$ will all have degree $O_{d,\ve}(1)$.
Before carrying out this plan, we
first show how (\ref{end-game}) will follow.  
Set $\#S_{p,q}(V;B,\pi,\om)=N$, and consider the matrix
$$
\mathbf{M}=\big(m_i(x_1^{(j)},x_2^{(j)},x_3^{(j)})\big)_{1\leq i \leq
  k,  ~1\leq j \leq N}.
$$
Then either $N<k$, or we may assume that $N \geq k$ and so $\D=0$ by
assumption. In either case it follows therefore that $\mathbf{M}$ 
has rank at most $k-1$.
Hence we may find a non-trivial linear combination $M$ of the monomials
$M_1,\ldots, M_k$ which vanishes at every point of
$S_{p,q}(V;B,\pi,\om)$, but not at all points of $Y$.  But then an application
of B\'ezout's theorem yields 
$$
\#S_{p,q}(V;B,\pi,\om)\leq
e\deg(M),
$$ 
which thereby establishes (\ref{end-game}).

It remains to establish that $\D=0$ provided that $k$ is chosen to be
sufficiently large in terms of $d$ and $\ve$.  For this we 
use the fact that (\ref{reduc}) holds for the points $x_1,\ldots,x_k$
that we are interested in, where
$\pi=[1,\pi_1,\pi_2,\pi_3]\in \bfP^3(\F_p)$ 
is a non-singular point on $X_p$ which is of multiplicity at most $2$
on the intersection with the tangent plane to $X_p$ at $\pi$, and 
$\om=[1,\om_1,\om_2,\om_3]\in \bfP^3(\F_q)$ is a non-singular point on
the curve $Y_q$.  The main obstacle that we shall have to overcome is
the fact that $\pi$ need not be a non-singular point on the curve $Y_p$,
despite being non-singular on $X_p$.  It is in the following
result that we make crucial use of the fact that 
$\pi$  is of multiplicity at most $2$
on the tangent plane section.

\begin{lem}\lab{lower}
For any $k\in \N$ we have $p^{\al(k)}q^{\be(k)}\mid \D$, where
$$
\al(k)\geq \frac{k^2}{2(e-1)} - O_{d,e,\ve}(k), \quad
\be(k) = \frac{k(k-1)}{2}.
$$
\end{lem}

\begin{proof}

Let $T_i=X_i/X_0$ for $1\leq i\leq 3$, and let $I\subset
\Q[T_1,T_2,T_3]$ be the ideal generated by the  
polynomials whose homogenisation belongs to the ideal $J$ of $Y$.
We begin by showing that $p^{\al(k)}\mid \D$.  Let $I_p$ be the image
of $I\cap\Z_{(p)}[T_1,T_2,T_3]$ in $\F_p[T_1,T_2,T_3]$, where $\Z_{(p)}$
denotes the localisation of $\Z$ at the prime $p$.
Let $A=A_p$ be the stalk of the curve $Y_p$ at the point
$\pi=[1,\pi_1,\pi_2,\pi_3] \in Y_p(\F_p)$, and let $\fr{m}$ be its
maximal ideal.

We now apply a rather general result due to the third author,
which can be found in the appendix.
Let $(n_\ell(A))_{\ell=1}^{\infty}$ be the sequence of
integers $n\geq 0$, in which each $n$ occurs precisely 
$\dim_{A/\fr{m}}\fr{m}^{n}/\fr{m}^{n+1}$ times.  
Then it follows that $p^{\mcal{A}(k)}$ divides $\D$, where
\begin{equation}\lab{a(k)}
\mcal{A}(k)=n_1(A)+\cdots+n_k(A).
\end{equation}
We need to show that 
\begin{equation}\lab{a(k)'}
\mcal{A}(k)\geq \frac{k^2}{2(e-1)} - O_{d,e,\ve}(k),
\end{equation}
for any $k\in \N$.  Let $B=B_p$ be the 
stalk of $X_p$ at $\pi$, and let $\fr{n}$ be the 
maximal ideal of this local ring.  Then $B$ is a regular local ring of
Krull dimension $2$ since $\pi$ is a non-singular point on $X_p$.
There is also a surjective ring homomorphism $\kappa: B\rightarrow A$ induced
by the embedding of $Y_p$ in $X_p$. Hence it follows that 
$$
\dim_{A/\fr{m}}\fr{m}^{n}/\fr{m}^{n+1}\leq \dim_{B/\fr{n}}\fr{n}^{n}/\fr{n}^{n+1}=n+1,
$$
for any $n\in \N$.
We now show that there exists an integer $n_0$ depending only on $d$
and $\ve$, such that
\begin{equation}\lab{arabian}
\dim_{A/\fr{m}}\fr{m}^{n}/\fr{m}^{n+1}\leq e-1,
\end{equation}
for all $n\geq n_0$.  It is clear that this will suffice to establish
(\ref{a(k)'}), via (\ref{a(k)}).

Recall that the prime ideal $J\subset \Q[X_0,X_1,X_2,X_3]$ of $Y$ is generated by $O_{d,\ve}(1)$ 
forms of maximal degree $O_{d,\ve}(1)$.  
We proceed by showing that the same is true of the ideal $J_p \subset \F_p[X_0,X_1,X_2,X_3]$. 
Now there are only $O_{d,\ve}(1)$ possible Hilbert functions for $J$
by \cite[Lemma 1]{broberg}, and it is not hard to see that 
$J \subset \Q[X_0,X_1,X_2,X_3]$ and 
$J_p \subset \F_p[X_0,X_1,X_2,X_3]$ have the same
Hilbert function since $(J\cap \Z_{(p)}[X_0,X_1,X_2,X_3])_D$ is a pure
sublattice of $(\Z_{(p)}[X_0,X_1,X_2,X_3])_D$ for each grade $D$.
It follows that there are only $O_{d,\ve}(1)$ possible Hilbert
functions for $J_p$, and so only $O_{d,\ve}(1)$ possible Hilbert
polynomials.  In order to deduce that $J_p$ is generated by forms of 
maximal degree $O_{d,\ve}(1)$, it is therefore enough to prove that 
$J_p$ is generated by forms of degree bounded solely in terms of
the Hilbert polynomial for $J_p$.  But this is a standard fact from the
theory of Hilbert schemes (see \cite[Theorem I.1.5]{kollar}, for example).

Putting everything together we therefore deduce that the kernel of
the homomorphism $\kappa: B\rightarrow A$ is generated by
$s=O_{d,\ve}(1)$ elements of $\fr{n}\setminus \fr{n}^{j+1}$ for 
some $j$ depending only on $d$.  Now let $A'$ (resp. $B'$) be the completion of
$A$ (resp. $B$) with respect to $\fr m$ (resp. $\fr n$), and let 
$\kappa': B'\rightarrow A'$ be the surjective homomorphism induced by
$\kappa$.   Then $B'$ is isomorphic to the ring $\F_p[[z_1,z_2]]$ of
formal power series, and the kernel of $\kappa'$ is generated by $s$
formal power series $\gamma_1(z_1,z_2), \ldots, \gamma_s(z_1,z_2)$
with a lowest term of degree at most $j$.  

We now show that $\dim_{A/\fr{m}}\fr{m}^{n}/\fr{m}^{n+1}$ is constant 
for all $n\geq n_0=2j$.  Our proof is based upon the observation that 
$\dim_{A/\fr{m}}\fr{m}^{n}/\fr{m}^{n+1}$ remains the same if we
replace $A$ by its completion $A'=\F_p[[z_1,z_2]]/(\gamma_1,\ldots,
\gamma_s)$. 
For each $1\leq j \leq s$ we may clearly replace $\gamma_j$ by any $\F_p$-linear
combination $\la_{1}\gamma_1+\cdots+\la_{j-1}\gamma_{j-1}+\gamma_j$
without changing the ideal $(\gamma_1,\ldots,\gamma_s)$.  We may
therefore assume that the smallest monomials in the power series
expansions of $\gamma_1,\ldots,\gamma_s$ are different.  
Here we employ the term ``smallest'' in the sense of the graded lexicographical ordering
$1,z_1,z_2,z_1^2,z_1z_2,z_2^2, \ldots$ of the monomials in $z_1$ and
$z_2$.  It follows that we may replace the set of generators 
$\gamma_1,\ldots,\gamma_s$ by the smallest monomials in their power series
expansions without changing the value of $\dim_{A/\fr{m}}\fr{m}^{n}/\fr{m}^{n+1}$.  
It therefore suffices to show that 
$\dim_{A/\fr{m}}\fr{m}^{n}/\fr{m}^{n+1}$ is constant for $n\geq 2j$
for a local ring $A=\F_p[[z_1,z_2]]/(\gamma_1,\ldots,
\gamma_s)$ of dimension $1$, where 
$\gamma_1,\ldots,\gamma_s$ are monomials in $z_1,z_2$ of degree at
most $j$.  Suppose now that $z_1^\mu$ and $z_2^\nu$ are the highest
powers of $z_1$ and $z_2$ that divide all of these monomials.
Consider the two sets
$$
\{z_1^n,z_1^{n-1}z_2, \ldots, z_1^{n-\nu+1}z_2^{\nu-1}\}, \quad
\{z_1^{\mu-1}z_2^{n-\mu+1},\ldots, z_1z_2^{n-1},z_2^n\},
$$
the first (resp. second) of which we define to be empty if $\nu=0$
(resp. $\mu=0$).
Then the union of these two sets 
form a basis for $\fr{m}^{n}/\fr{m}^{n+1}$ whenever $n\geq 2j$.
Hence $\dim_{A/\fr{m}}\fr{m}^{n}/\fr{m}^{n+1}=\mu+\nu$ for $n\geq 2j$.

We may now complete the proof of (\ref{arabian}).  Since $Y$ has
dimension $1$ it follows that the multiplicity $\rom{mult}_{\pi}(Y_p)$
of $Y_p$ at $\pi$ is equal to the dimensions of the vector spaces
$\fr{m}^{n}/\fr{m}^{n+1}$ over $A/\fr{m}$ for sufficiently large $n$
(see Fulton \cite[Excercise 4.3.1]{fulton}, for example).  
We clearly have the trivial upper bound 
$$
\rom{mult}_{\pi}(Y_p)\leq \deg Y_p=\deg Y=e.
$$  
Suppose now that  $\rom{mult}_{\pi}(Y_p)=e.$  Then $Y_p$ must be a
cone with vertex $\pi$.  In particular $Y_p$ is contained in the
intersection of $X_p$ with its tangent plane at $\pi$.  This
contradicts the assumption that $\pi$ has multiplicity at most $2<e$
on the tangent plane section of $X_p$ at $\pi$, and so establishes (\ref{arabian}).

It remains to show that $q^{\be(k)}\mid \D$.  
Indeed, since $p\not\in \Omega$ by Lemma \ref{q} 
it will then follow that 
$p^{\al(k)}q^{\be(k)}\mid \D$, as required.
Adopting the same notation as above, we let 
$A=A_q$ be the stalk of the curve $Y_q$ at the non-singular point
$\om=[1,\om_1,\om_2,\om_3] \in Y_q(\F_q)$, and let $\fr{m}$ be its
maximal ideal.
Then an application of the result in the appendix yields
$$
q^{n_1(A)+\cdots+n_k(A)}\mid \D.
$$  
It is then a simple matter to note that 
$n_1(A)+\cdots+n_k(A)=\be(k)$ in this case.  Indeed 
for any $n\geq 0$ the dimension of $\fr{m}^{n}/\fr{m}^{n+1}$ over
$A/\fr{m}$ is now equal to the
multiplicity  of the non-singular point $\om$ on the reduced curve $Y_q$,
which is plainly equal to $1$.
\end{proof}

We may now complete the proof of Theorem \ref{ash-2}.  Recall the
inequality (\ref{prime}) for $p$.  
On combining
Lemma \ref{lower} with the upper bound (\ref{D1}) for $\log |\D|$, we
easily deduce that $\D=0$ whenever 
$$
\log q >
\left(\frac{1}{e}-
  \frac{1}{(e-1)\sqrt{d}}\right)\log B -\frac{\ve \log B}{e-1}+  O_{d,e,\ve}(1)
+O_{d,e,\ve}\left(\frac{\log B}{k-1}\right).
$$
Upon taking $k$ to be sufficiently large compared with $d$ and $\ve$, 
this is plainly satisfied for every $q\in \Omega$, by (\ref{q-prime}).

\appendix

\section{Appendix}

Let $p$ be a prime and let $\Z_{(p)}$ denote the localisation of $\Z$ at
$p$. Let $R$ be any commutative noetherian local ring containing
$\Z_{(p)}$ as a subring.
We write $A$ for the quotient ring $R/pR$, and let $\fr{m}$ be the
maximal ideal of $A$.
Finally, let 
$$
n_1(A),n_2(A), n_3(A),\ldots
$$ 
be the non-decreasing
sequence of integers $n\geq 0$  in which each $n$ occurs precisely 
$\dim_{A/\fr{m}}\fr{m}^{n}/\fr{m}^{n+1}$ times.  
The purpose of this appendix is to record the following result due to
the third author \cite{s}.

\begin{lemma}
Let $\psi_1,\ldots,\psi_k$ be ring homomorphisms from $R$ to $\Z_{(p)}$,
for some positive integer $k$, and let $r_1,\ldots, r_k \in R$.
Then it follows that 
$$
p^{n_1(A)+\cdots+n_k(A)} \mid \det \big(\psi_i(r_j)\big)_{1\leq i,j
  \leq k}.
$$
\end{lemma}

The main idea in the proof of this result is to replace each $r_j$ by
a suitable linear combination $\la_1r_1+\cdots+\la_{j-1}r_{j-1}+r_j$,
with coefficients in $\Z_{(p)}$, such that $\psi_1(r_j), \ldots,
\psi_k(r_j)$ are all divisible by $p^{n_j(A)}$.

\end{document}